\documentclass{amsart}
\usepackage{amsopn,amsfonts,amscd,amssymb,latexsym,amsmath,amsthm}

\newtheorem{theorem}{Theorem}[section]
\newtheorem{remark}{Remark}[section]
\newtheorem{example}{Example}[section]
\newtheorem{lemma}{Lemma}[section]
\newtheorem{problem}{Problem}[section]
\newtheorem{corollary}{Corollary}[section]
\newtheorem{definition}{Definition}[section]
\newtheorem{condition}{Condition}[section]

\renewcommand{\i}[1]{\raisebox{1pt}{$\stackrel{\scriptscriptstyle \circ}{#1}$}}

\begin{document}

\title[On the Cauchy  Problem  for the  Elliptic Complexes]
{On the Cauchy  Problem  
for the  Elliptic Complexes  in Spaces of Distributions}

\author[A. Shlapunov]{Alexander Shlapunov}
\address{Siberian Federal University,
         Institute of Mathematics and Computer Science,
         pr. Svobodnyi 79,
         660041 Krasnoyarsk,
         Russia}
\email{ashlapunov@sfu-kras.ru}

\author[D. Fedchenko]{Dmitrii Fedchenko}
\address{Siberian Federal University,
         Institute of Core Undergraduate Programs,
         pr. Svobodnyi 79,
         660041 Krasnoyarsk,
         Russia}
\email{dfedchenk@gmail.com}

\begin{abstract}
Let $D$ be a bounded domain in ${\mathbb R}^n$, $n\geq 2$,  with a smooth 
boundary $\partial D$.  We indicate appropriate Sobolev spaces of negative 
smoothness to study the non-homogeneous Cauchy  problem  for an  elliptic 
differential complex $\{A_i\}$ of first order operators. In particular, we 
describe traces on $\partial D$ of tangential part $\tau_i (u)$ and normal part
$\nu_i(u)$ of a (vector)-function $u$ from the corresponding Sobolev space and 
give an adequate formulation of the problem. If the Laplacians of the complex 
satisfy the uniqueness condition in the small then we obtain necessary and 
sufficient solvability conditions of the problem and produce formulae for its 
exact and approximate solutions.  For the Cauchy problem in 
the Lebesgue space $L^2(D)$ we construct the approximate and exact solutions 
to the Cauchy problem with the maximal possible regularity. Moreover, using 
Hilbert space methods, we construct Carleman's  formulae for 
a (vector-) function $u$ from the Sobolev space $H^{1}(D)$ by its Cauchy data 
$\tau _i (u)$  on a subset $\Gamma \subset \partial D$ and the values of $A_i 
u$ in $D$  modulo the null-space of the Cauchy problem. Some instructive 
examples for elliptic complexes of operators with constant coefficients are 
considered.
\end{abstract}

\maketitle
\medskip

\noindent \underbar{Key words:} Elliptic differential complexes, ill-posed  
Cauchy problem, Carleman's  formula.

\medskip
 
It is well known that the Cauchy problem for an elliptic system $A$ is 
ill-posed (see, for instance, \cite{Hd}). Apparently, the serious investigation of 
the problem was stimulated by practical needs. Namely, it naturally appears in 
applications: in hydrodynamics (as the Cauchy problem for holomorphic 
functions), in geophysics (as the Cauchy problem for the Laplace operator), in 
elasticity theory (as the Cauchy problem for the Lam\'e system) etc., see, for 
instance, the book \cite{Tark37} and its bibliography. The problem was 
actively studied through the XX century (see, for instance, \cite{La}, 
\cite{Maz}, \cite{KoLan},  \cite{A},  \cite{Na}, \cite{AKy}, \cite{ShTaLMS}, 
\cite{SchShTa2},  \cite{ShTaRJMP}, \cite{SheShlL2}, \cite{SheShl3} and many others). 

Differential complexes appear as compatibility conditions for overdetermined 
operators (see, for instance, \cite{Spe}, \cite{Ta1}). Thus, the Cauchy 
problem for them is of the special interest. One of the first problems of this 
kind was the Cauchy problem for the Dolbeault complex (the compatibility 
complex for the multidimensional Cauchy-Riemann system), see \cite{AnHi}. The 
interest to it was great because of the famous example by H.~Lewy of the 
differential equation without solutions, constructed with the use of the 
tangential Cauchy-Riemann operator, see \cite{Lewy}. Recently new approaches 
to the problem were found in spaces of smooth functions (see \cite{NaSchTa}, 
\cite{BrinkHill}) and in spaces of distributions (see \cite{FeShl1}, \cite{FeShl2}). 

We consider the Cauchy problem in spaces of distributions with some restrictions on growth in 
order to correctly define its traces on boundaries of domains (see, for instance,  
\cite{Tark37}, \cite{Sche2}, \cite{Stra}, \cite{Roj}, \cite{LiMa},). 
In this paper we develop the approach 
presented in \cite{ShTaLMS} to study the {\it homogeneous} Cauchy problem for 
overdetermined elliptic partial differential operators. Instead we consider the 
{\it non-homogeneous} Cauchy problem for elliptic complexes. 

\section{Preliminaries} \label{s.1}

\subsection{Differential complexes}

Let $X$ be a $C^{\infty}$-manifold of dimension $n\geq 2$ with a smooth 
boundary $\partial X$. We tacitly assume that it is enclosed into a smooth 
closed manifold $\tilde{X}$ of the same dimension.

For any smooth $\mathbb C$-vector bundles $E$ and $F$ over $X$, we write 
$\mbox{Diff}_{m} (X; E \to F)$ for the space of all the linear partial 
differential operators of order $\leq m$ between sections of $E$ and $F$. 
Then, for an open set $O\subset \i{X}$ (here $\i{X}$ is the interior of $X$)  
over which the bundles and the manifold are trivial, the sections over $O$ may 
be interpreted as (vector-) functions and $A\in \mbox{Diff}_m(X; E\rightarrow 
F)$ is given as $(l\times k)$-matrix of scalar differential operators, i.e. we 
have 
$$
A=\sum_{|\alpha|\leq m}a_\alpha(x) \frac{\partial^{|\alpha|}}{\partial 
x^\alpha}, \quad x \in O, 
$$ 
where $a_\alpha(x)$ are $(l\times k)$-matrices of 
$C^\infty (O)$-functions, $k =\mbox{rank}(E)$, $l=\mbox{rank}(F)$.

Denote $E^{\ast}$ the conjugate bundle of $E$.
Any Hermitian metric $(.,.)_{x}$ on $E$ gives rise to a sesquilinear bundle
isomorphism (the Hodge operator) $\star_{E} \! : E \to E^{\ast}$ by the 
equality $\langle \star_{E} v, u \rangle_{x} = (u,v)_{x}$ for all sections $u$ 
and $v$ of $E$; here $\langle ., . \rangle_{x}$ is the natural pairing in the 
fibers of $E^*$ and $E$. 

Pick a volume form $dx$ on $X$, thus identifying the dual and the conjugate bundles.
For $A \in \mbox{Diff}_{m} (X; E \to F)$, denote 
by $A^{\ast} \in \mbox{Diff}_{m} (X; F \to E)$ the formal adjoint operator. 

Let $\pi: T^{\ast} X \to X$ be the (real) cotangent bundle of $X$ and 
let $\pi^*E $ be a induced bundle for the bundle $E$ (i.e. the fiber 
of $\pi^*E$ over the point $(x,z) \in T^{\ast} X $ coincides with $E_x$). 
We write $\sigma  (A): \pi^*E  \to \pi^*F $ for the principal homogeneous 
symbol of the order $m$ of the operator $A$. 
 
Let $D$ be a bounded domain (i.e. open connected set) in $\i{X}$ with 
infinitely differentiable boundary $\partial D$. Denote $C^\infty (D,E)$ 
the Fr\'echet space of all the infinitely differentiable sections of the 
bundle $E$ over $D$ and denote $C^\infty  (\overline D,E)$ the subset in 
$C^\infty (D,E)$ which consists of sections with all the derivatives 
continuously extending up to $\overline D$. Let also $C^\infty _{comp} (D,E)$ 
stand for the set of all the smooth sections with compact supports in $D$. 
Besides, for an open (in the topology of $\partial D$) subset  $\Gamma \subset 
\partial D$, let $C^\infty_{comp} (D \cup \Gamma,E)$ be the set of  all the 
$C^\infty(\overline D,E)$-sections with compact supports in $D \cup \Gamma$. 

For a distribution-section $u \in (C^{\infty}_{comp} (D, E))'$ we always 
understand $Au$ in the sense of  distributions in $D$. The spaces of all the 
weak solutions of the operator $A$ in $D$ we denote $S_A(D)$. 

We often refer to the so-called uniqueness condition in the small on $\i{X}$ 
for an operator $A$. 

\begin{condition} \label{US}
If $u$ is a distribution in a domain $D\Subset \i{X}$ with $Au =0$ in $D$ 
and $u = 0$ on an open subset $O$ of  $D$ then $u \equiv 0$  in $D$.
\end{condition}
It holds true if, for instance, all the objects under consideration are real 
analytic.

Let $G_A(.,.)\in \mbox{Diff}_{m-1}(X;(F^*,E)\rightarrow\Lambda^{n-1})$ denote 
a Green operator attached to $A$, i.e. such a bi-differential operator that 
$$
dG_A(\star_F g,v) = ((Av,g)_x - (v, A^* g)_x ) \,dx
\mbox{ for all } g\in C^\infty(X,F), \ v\in C^\infty(X,E);
$$
here 
$\Lambda^{p}$ is the bundle of the exterior differential forms of the degree 
$0\leq p \leq n$ over $X$. The Green operator always exists (see \cite
[Proposition 2.4.4]{Ta1}) and for the first order operator $A$  it may be 
locally written in the following form:
$$
G_A(\star g,v) = g^*(x) \ \sigma (A)(x, (\star dx_1, \dots, \star dx_n)) \ v 
(x) \mbox{ for all }  g\in C^\infty(X,F), \ v\in C^\infty(X,E).
$$
Then it follows from Stokes formula that 
the (first) Green formula holds true:
\begin{equation} \label{eq.Green1}
\int_{\partial D} G_A(\star g,v) 
= \int_D ( (Av,g)_x - (v, A^*g)_x ) \ dx 
\mbox{ for all } g\in C^\infty(X,F), \ v\in C^\infty(X,E).
\end{equation}

Fix a defining function of the domain $D$, i.e. a real valued $C^\infty$-smooth 
function $\rho$  with $|\nabla \rho| \ne 0$  on $\partial D$ and such that
$D = \{ x \in X: \ \rho(x) < 0 \}$. Without loss of a generality
we can always choose the function $\rho$ in such a way that $|\nabla
\rho| = 1$ on a neighborhood of $\partial D$.
Then 
\begin{equation} \label{eq.Green.op}
G_A(\star g,v) = \int_{\partial D}  (\sigma (A)(x, \nabla \rho) 
\ v , g)_x  \  \ ds (x) 
\mbox{ for all }  g\in C^\infty(X,F), \ v\in C^\infty(X,E), 
\end{equation}
where  $ds$ is the volume form on $\partial D$ induced from $X$. 

Our principal object to study will be a complex $\{A_i,E_i\}_{i=0}^N$ of partial 
differential operators over $X$ (see, \cite{Ta1}, \cite{Spe}), 
\begin{equation} \label{eq.complex}
0\rightarrow 
C^\infty(X,E_0)\stackrel{A_0}{\rightarrow}C^\infty(X, E_1)\stackrel{A_1}{
\rightarrow}C^\infty(X, E_2) \rightarrow \dots 
\stackrel{A_{N-1}}{\rightarrow} C^\infty(X, E_N) \rightarrow 0, 
\end{equation}
where $E_i$ are the bundles over $X$ and $A_i \in \mbox{Diff}_{1} (X; E_i \to 
E_{i+1})$ with $A_{i+1} \circ A_i \equiv 0$; we tacitly assume that 
$A_i =0$ for both $i<0$ and $i\geq N$. Obviously, $\sigma(A_{i+1}) \circ 
\sigma(A_i) \equiv 0$. We say that the complex 
$\{A_i,E_i\}_{i=0}^N$  is elliptic if the corresponding symbolic complex, 
\begin{equation} \label{eq.complex.symb}
0\rightarrow 
\pi^* E_0\stackrel{\sigma(A_0)}{\rightarrow} \pi^* E_1\stackrel{ \sigma 
(A_1)} {\rightarrow} \pi^* E_2 \rightarrow \dots  
\stackrel{ \sigma (A_{N-1})} {\rightarrow}  \pi^* E_{N} \rightarrow 0,
\end{equation}
is exact for all $(x,z) \in T^{\ast} X \setminus \{0\} $, i.e. the the range of the map 
$\sigma (A_i)$ coincides with the kernel of the map $\sigma (A_{i+1})$. 
In particular, $\sigma (A_0)$ is injective away from the zero
section of $T^{\ast} X$ and $\sigma (A_{N-1})$ is surjective. 

As any differential complex is homotopically equivalent to 
a first order complex, we will consider elliptic complexes of first order 
operators only. Hence it follows that the Laplacians $\Delta_i = 
A_i^{*} A_i + A_{i-1} A_{i-1}^{*} $ of the complex are elliptic differential 
operators of the second order and types $E_i \to E_i$ on $X$ for 
$0\leq i \leq N$. 

\subsection{Sobolev spaces}

We write $ L^{2} (D, E)$ for the Hilbert space of all the measurable sections 
of $E$ over $D$ with a scalar product $(u, v)_{L^{2} (D, E)} = \int_{D} 
(u,v)_{x} dx$. We also denote $H^{s} (D, E)$ the Sobolev space
of the distribution sections of $E$ over $D$, whose weak derivatives up to the 
order $s \in {\mathbb N}$ belong to $L^{2} (D, E)$. As usual,  let $H^s 
_{loc}(D\cup \Gamma,E)$ be the set of sections in $D$ belonging to $H^s 
(\sigma,E)$ for every measurable set $\sigma$ in $D$ with $\overline \sigma 
\subset D\cup \Gamma$. 

Further, for non-integer positive $s$ we define the Sobolev spaces $H^{s} (D, E)$ 
with the use of the proper interpolation procedure (see, for example, \cite[\S 
1.4.11]{Tark37}). In the local situation we can use other (equivalent)
approaches. For instance, if $X \subset {\mathbb R}^n $ and the bundle $E$ is 
trivial, we may  we denote $H^{1/2} (D,E)$ the closure of $C^\infty (\overline 
D,E)$ functions with respect to the norm (see \cite{EgShb}):
$$
\| u \|_{H^{1/2} (D,E)} = \sqrt{ \| u \|^2 _{L^{2} (D,E)} +
\int_{D}   \int_{D}  \frac{|u(x)- u(y)|^2 d x \, dy }{|x-y|^{2n+1}}}.
$$ 
Then, for $s \in {\mathbb N}$,  let 
$H^{s-1/2} (D,E)$ be the space of functions from $H^{s-1} (D,E)$ such that 
the weak derivatives of the order $(s-1)$ belong to $H^{1/2} (D, E)$.

The Sobolev spaces of negative smoothness are usually defined with the use 
of a proper duality (see \cite{Ad}). For instance, one can 
consider the Sobolev space  $\tilde H^{-s} (D,E)$ as the completion of 
the space $C^{\infty} _{comp}(D,E)$ with respect to the norm 
$\sup\limits_{v \in C^{\infty} _{comp} (D,E)} \frac{|(u,v)_{L^{2} (D,E)}|}{\| 
v \|_{H^{s}(D,E)}}$, $s\in \mathbb N$. Unfortunately, elements of these spaces 
may have "bad"{} behavior near $\partial D$, but the study of the Cauchy 
problem needs a correctly defined notion of a trace. This is the reason 
we use slightly different spaces; we follow \cite{Sche2} (cf. \cite{Roj}, 
\cite[Chapters 1, 9]{Tark37}, \cite{ShT5}). More exactly, denote by 
$C^{\infty}_{m-1} (\overline D,E)$ the subspace in $C^{\infty} 
(\overline D,E)$ consisting of the sections with vanishing on $\partial D$
derivatives up to order $m-1$. Let $s \in {\mathbb N}$. For sections  
$u \in C^{\infty} (\overline D,E)$ we define two types of negative norms 
$$
{}\| u \|_{-s}  = \sup_{v \in C^{\infty} (\overline D,E)}
\frac{|(u,v)_{L^{2} (D,E)}|}{\| v \|_{H^{s}(D,E)}}, \quad 
{}| u |_{-s}  = \sup_{v \in C^{\infty}_{m-1}(\overline D,E)}
\frac{|(u,v)_{L^{2} (D,E)}|}{\| v \|_{H^{s}(D,E)}}.
$$
It is more correct to write $\| \cdot \|_{-s,D,E}$ and $| \cdot |_{-s,D,E}$, but 
we prefer to omit the indexes $D$, $E$, if it does not cause misunderstandings.  
It is convenient to set $\| \cdot \|_{0,D}= \| \cdot \|_{L^2(D,E)}$.

Denote the completions of space $C^{\infty} (\overline D,E)$ with respect to 
these norms by $H^{-s} (D,E)$ and $H (D,E, | \cdot |_{-s})$ respectively. It 
follows from the definition that the elements of these Banach spaces are 
distributions of finite orders on $D$ and these spaces could be called the Sobolev 
spaces of negative smoothness. Clearly, they satisfy the following relations:  
$H^{-s} (D,E)\hookrightarrow H (D,E, | \cdot |_{-s})
\hookrightarrow \tilde H^{-s} (D,E)$, and, similarly, 
$H^{-s} (D,E)\hookrightarrow H^{-s-1} (D,E)$, 
$H (D,E, | \cdot |_{-s})\hookrightarrow H (D,E, | \cdot |_{-s-1})$.

The Banach space  $H^{-s} (D,E)$  
can be identifyed with the dual space $\left( 
H^{s} ({D},E) \right)'$ of the standard Hilbert space $H^{s} ({D},E)$ (see,  
for instance, \cite[Theorem 1.4.28]{Tark37}). 

Clearly, any element $u\in H ^{-s}(D,E)$ extends up to an element 
$U \in H ^{-s}(\i{X}, E)$ via 
$$ 
\langle U, v\rangle_{\i{X}} =\langle u, v\rangle _{D}
\mbox{ for all } v \in H^{s} (\i{X},E); 
$$
here $\langle \cdot, \cdot \rangle_D$ is a pairing $H \times H'$ for  
a space $H$ of distributions over $D$. It is natural to denote this extension 
$\chi_{D} u$ because its support belongs to $\overline D$. Obviously, this 
extension induces a bounded linear operator 
\begin{equation} \label{eq.chi}
\chi_{D}: H ^{-s}(D,E) \to H ^{-s}(\i{X},E), \quad s  \in {\mathbb Z}_+.
\end{equation}

It is known that the differential operator $A$ continuously maps $H^{s} (D,E)$ 
to $H^{s-m} (D,F)$, $m \leq s$, $s \in \mathbb N$. The following lemma shows
the specific way of the action of $A$ for $s\leq 0$.

\begin{lemma} \label{l.A} A differential operator $A$ induces linear
bounded operator $A: H ^{-s} (D,E) \to H (D,F, | 
\cdot |_{-s-m})$, $s  \in {\mathbb Z}_+$. 
\end{lemma}

{\bf Proof.} Immediately follows from (\ref{eq.Green1}) and 
(\ref{eq.Green.op}). 
\hfill $\square$

However there is no need for elements of $H ^{-s}(D, E)$ to have a trace on 
$\partial D$ and there is no  need for $A$ to map  $H ^{-s}(D,E)$ to 
$H ^{-s-m}(D, F)$. 

\section{Traces of Sobolev functions of negative smoothness}

By the discussion above we need to introduce some other spaces in  order to 
define the traces on $\partial D$. In general, our approach is closed to 
the one described in \cite[\S 9.2, 9.3]{Tark37}.

\subsection{Strong traces on the boundary}
 
It is well known that if $\partial D$ is sufficiently smooth then 
the functions from the Sobolev space $H^s (D)$, $s \in \mathbb N$, have traces 
on the boundary in the Sobolev space $H^{s-1/2} (\partial D)$ and the 
corresponding trace operator $t_s: H^s(D) \to H^{s-1/2} (\partial D)$ is 
bounded and it admits the bounded right inverse operator (see, for instance, \cite{EgShb}).
 In particular, 
this means that for every $u \in H^s_{loc} (D\cup \Gamma,E)$, 
$s \in \mathbb N$, there is the trace $t_{\Gamma,E}(u)$ on $\Gamma$ belonging to 
$H^{s-1/2}_{loc} (\Gamma,E)$.

In order to define the so-called strong traces on $\partial D$ for elements of 
Sobolev spaces with negative smoothness we denote $H ^{-s}_{t}(D,E)$ the 
completion of $C^{\infty} (\overline D,E)$ with respect to the graph-norm:
$$ 
\| u \|_{-s,t} = \left( \| u \|^2_{-s}+  \| u \|^2_{-s-1/2,\partial D}\right) 
^{1/2}.
$$

Thus the operator $t_s$ induces the bounded linear trace operator  
$$ 
t_{-s,E}: H^{-s}_{t} (D, E) \to 
H ^{-s-1/2}(\partial D, E). 
$$

\begin{remark} \label{r.Dirichlet}
{\rm The spaces $H^{-s} (D,E)$,  $H^{-s}_t (D,E)$, 
$H (D,E, | \cdot |_{-s})$ are well known. Let $A$ be a first order operator 
with injective principal symbol. Given distributions $w$ and $u_0 $, consider 
the Dirichlet problem for strongly elliptic formally self-adjoint  second 
order operator $A^*A$. It consists in finding a distribution $u$ satisfying  
\begin{equation} \label{eq.HDP}
\left\{ 
\begin{array}{rclcl}
A^*A u & = & w & \mbox{in} & D,\\
t (u)    & = & u_0 & \mbox{on} & \partial D.
\end{array} 
\right.                
\end{equation}
It follows from \cite[theorems 2.1 and 2.2]{Sche2} (see also \cite{Roj}, 
\cite{ShT5} for systems of equations) that the Uniqueness Theorem and the 
Existence Theorem are valid for problem (\ref{eq.HDP}) on the Sobolev 
scale $H^s (D,E)$, $s\in\mathbb Z$ for data $w \in H (D, E, | \cdot 
|_{s-2}) $ and $u_0 \in H^{s-1/2}(\partial D,E)$. Denote by ${\mathcal 
P}^{(D)}$ the operator mapping $u_0$ and $w=0$ to the unique solution to the 
Dirichlet problem (\ref{eq.HDP}). Similarly, denote ${\mathcal G} ^{(D)}$ 
the operator mapping $w$ to the unique solution to the Dirichlet problem 
(\ref{eq.HDP}) with zero boundary Dirichlet data. Clearly, ${\mathcal 
G} ^{(D)}_{A^*A}$ is the famous Green function of the Dirichlet problem (\ref
{eq.HDP}) and ${\mathcal P}^{(D)}_{A^*A}$ is the Poisson integral 
corresponding to the problem. The standard theorem of improving the smoothness 
of the Dirichlet problem (see,  for instance, \cite{EgShb} or \cite[Theorem 
9.3.17]{Tark37}) and \cite[Theorem 2.26 and Corollary 2.31]{ShT5} imply that 
the operators ${\mathcal P}^{(D)}_{A^*A}$, ${\mathcal G} ^{(D)}_{A^*A}$ act 
continuously on the following Sobolev scale: 
$$
{\mathcal P}^{(D,s)}_{A^*A}: H ^{s-1/2}(\partial D,E) \to H^s (D,E), \quad 
{\mathcal G} ^{(D,s)}_{s,A^*A}: H (D, E, | \cdot |_{s-2}) \to H^s (D,E), \quad 
s \in {\mathbb N},
$$
$$
{\mathcal P}^{(D,-s)}_{A^*A}: H ^{-s-1/2}(\partial D,E) \to H_t^{-s} (D,E), 
\quad {\mathcal G}^{(D,-s)}_{A^*A}: H (D, E, | \cdot |_{-s-2}) \to H_t^{-s} 
(D,E), \quad s \in {\mathbb Z}_+.
$$
They completely describe the solutions of the Dirichlet 
problem on the scale of the Sobolev spaces.
}
\hfill $\square$
\end{remark}

However we need a more subtle characteristic of the traces to study the Cauchy 
problem for the differential complex $\{A_i\}$. 

For a section $u$ of $E$ over $D$ and a first order operator $A$, let $\tilde 
\tau_A(u) = \sigma (A) (x, \nabla \rho (x)) u $ represent the Cauchy data of $u$
with respect to $A$ (see, for instance, \cite[\S 3.2.2]{Ta1}). Similarly, 
let $\tilde \nu _A(f) =  \tilde \tau_{A^*}(f) $ represent the 
Cauchy data of $f$ with respect to $A^*$ for a section $f$ of $F$. Then the 
maps  $\tilde \tau$, $\tilde \nu $ induces a bounded linear operators 
\begin{equation} \label{eq.tau.s}
\tilde \tau_{A,s} : H^s (D,E)  \to  H^{s-1/2} (\partial D, F), 
\quad \tilde \nu_{A,s} : H^s (D,F)  \to  
H^{s-1/2} (\partial D, E), \quad s \in {\mathbb N}.  
\end{equation} 

Denote the completions of the space $C^{\infty} (\overline D, E)$
with respect to graph-norms
$$ 
\| u \|_{-s,A} = \left( \| u \|^2_{-s} + \| A u \|^2_{-s-1}\right)  ^{1/2}, 
\quad  \| u \|_{-s, \tilde \tau_A} = \left( \| u \|^2_{-s}+ \|\tilde \tau _A(u) 
\|^2_{-s-1/2,\partial D}\right) ^{1/2}
$$
by 
$H^{-s}_{A} (D,E)$ and $H ^{-s}_{\tilde \tau_A}(D,E)$ respectively.
Clearly, the elements of these spaces are more regular in $\overline D$ than 
the elements of $H ^{-s} (D,E)$. Moreover, by the very definition, the differential 
operator $A$ induces a bounded linear operator 
$$
A _{-s}: H^{-s}_{A} (D, E) \to 
H ^{-s-1}(D,F), 
$$
and the trace operator (\ref{eq.tau.s}) induces a bounded linear operator   
$$ 
\tilde \tau _{A,s}: H^{-s}_{\tilde \tau_A} (D,E) \to 
H^{-s-1/2}(\partial D, F). 
$$

\begin{theorem} \label{t.equiv.strong.tau}  
The linear spaces $H^{-s} _{A}(D, E)$ and $H ^{-s}_{\tilde \tau _A}(D,E)$ 
coincide and their norms are equivalent. Moreover, if $A$ has the injective 
principal symbol then the spaces $H^{-s} _{\tilde \tau_A}(D, E)$ and $H ^{-s}_{t}(D,E)$
coincide and their norms are equivalent.
\end{theorem}

{\bf Proof.} It follows from the definition of the spaces that we need to check 
the relations between the norms on the sections from $C^\infty (\overline D,E)$ only. 
By Green's formula (\ref{eq.Green1}) and (\ref{eq.Green.op}) we have for all 
$u \in C^\infty (\overline D, E)$: 
\begin{equation*} 
\|u  \|_{-s, A}^2  \leq ( 1+\| A_{s+1} ^*\|^2+ \| t_{s+1,F} \|^2) (\|u \|^2_
{-s} + \| \tilde \tau_A(u) \|_{-s-1/2,\partial D} ^2), 
\end{equation*}
where $A ^*_{s+1}  : H^{s+1} (D, F) \to  
H^{s}(D, F)$ is the linear bounded operator induced by 
the differential operator $A^*$. 

Back, fix a section $g_0 \in C^\infty (\partial D,F)$. 
Now let $\nabla_F \in \mbox{Diff}_1 (X; F\to F\otimes (T^*X)_c)$ 
be a connections in the bundle $F$ compatible with the corresponding Hermitian 
metric (see \cite[Ch. III, Proposition 1.11]{We}). Obviously $\nabla_F $
has the injective symbol. Then, using remark \ref{r.Dirichlet} we see that 
there is a section $g \in C^\infty (\overline D,F)$ with $g=g_0$ on $\partial D$ 
and  $\|g\|_{s+1} \leq 
\gamma \|g_0\|_{s+1/2}$. For instance we may take $g = {\mathcal P}
_{\nabla_F ^*\nabla_F}^{(D)} g_0$. 
Therefore  Green's formula (\ref{eq.Green1}) and formula (\ref{eq.Green.op}) 
imply that for all $u \in C^\infty (\overline D, E)$ we have:
$$ 
\int_{\partial D} ( \tilde \tau_A (u),  g_0)_x  ds (x) = 
\int_D ((Au, g)_x - (u,A^*g)_x) dx.
$$
Hence
\begin{equation*} 
\|u  \|_{-s}^2+ \|\tilde \tau_A (u) \|^2_{-s-1/2, \partial D} \leq(1+ \gamma ^2
\|A ^*_{s+1}\| ^2 +\gamma ^2 ) (\|u \|^2_{-s} +   \| A u \|^2_{-s-1} ), 
\end{equation*}
i.e. the spaces $H^{-s} _{A}(D, E)$ and $H ^{-s}_{\tilde \tau _A}(D,E)$ 
coincide and their norms are equivalent.

Finally, if the symbol $\sigma(A)$ is injective then the map 
$ \sigma^*(A) (x, \nabla \rho (x))\sigma(A)(x, \nabla \rho (x))$ 
is invertible on $\partial D$ and 
$$ 
\tilde \tau _A(u) = \sigma(A)(x, \nabla \rho (x)) t (u), \quad 
t (u) = (\sigma^*(A) (x, \nabla \rho (x))\sigma(A)(x, \nabla \rho (x)))^{-1} 
\tilde \nu _A (\tilde \tau _A(u)), 
$$
which means that the norms $\|\cdot\|_{-s,t}$ and $\|\cdot\|_{-s,\tilde \tau 
_A}$ are equivalent on $C^\infty (\overline D,E)$.
\hfill $\square$

Now for the complex $\{A_i\}$  denote $\tilde \tau_i$ the Cauchy data with 
respect to $A_i$. Similarly denote $\tilde \nu_i$ the Cauchy data with respect 
to $A^*_{i-1}$. As the complex is elliptic then the matrix 
$L(x)= \sigma^*(A_i) (x, \nabla \rho (x))\sigma(A_i)(x, \nabla \rho (x))
+ \sigma(A_{i-1}) (x, \nabla \rho (x))\sigma^* (A_{i-1})(x, \nabla \rho (x)) $ 
is invertible in a neighborhood of $\partial D$. Then we set 
$$
\tau_i = L^{-1}(x) \tilde  \nu_{i+1} \circ \tilde \tau_i, \quad  
\nu _i  = L^{-1}(x)  \tilde\tau_{i-1} \circ \tilde \nu_i .
$$

\begin{lemma} \label{l.projections}
The following identities hold true:
$$
\tilde \tau_{i+1} \circ \tilde \tau_{i} =0, \, \tilde \nu_{i-1} \circ \tilde 
\nu_{i} =0, \, \tilde \tau_{i} \circ \nu_{i} =0, \, \tilde \nu _{i} 
\circ \tau _i =0, \, \tilde \tau_i =  \tilde \tau_i \circ \tau_i,  \, 
\tilde \nu _i = \tilde  \nu _i\circ \nu_i , 
$$
$$
\tau_i  \circ  \tau _i = \tau_i, \, \nu _i \circ  \nu_i = \nu _i, \,  \tau_i 
\circ  \nu_i   =0, \,  \nu _i \circ  \tau_i  =0, \, \tau _i + \nu_i  = 1,
$$
$$
\tau_i ^*  = \tau_i, \, \nu _i ^*  = \nu _i, \,  \tilde \tau_i^* =   
 \tilde \nu_{i+1}  , \,  \tilde \nu _i ^* = \tilde \tau_{i-1} .
$$
\end{lemma}

{\bf Proof.} See, for instance, \cite[formulae (3.2.3)]{Ta1}. \hfill $\square$

Because of Lemma \ref{l.projections}, the projections $\tau_i (u)$ and $\nu_i 
(u)$ are often called the tangential and the normal parts of a section $u$ with 
respect to the complex $\{A_i\}$ respectively. 

Due to  Lemma \ref{l.projections} we have for all $ u\in C^\infty 
(\overline D, E_i)$,  $g\in C^\infty (\overline D, E_{i+1})$:
\begin{equation} \label{eq.Green.i}
\int_{\partial D} ( \tau_i (u),  \tilde \nu _{i+1} (g))_x  \ ds (x) = \int_D 
((A_iu, g)_x - (u,A^*_ig)_x) dx . 
\end{equation}
Denote the completion of the space $C^{\infty} (\overline D, E_i)$
($0\leq i \leq N$) with respect to graph-norms
$$ 
\| u \|_{-s, \tau_i} = \left( \| u \|^2_{-s}+ \|\tau_i (u) 
\|^2_{-s-1/2,\partial D}\right) ^{1/2}, \quad
\| u \|_{-s, \nu_i} = \left( \| u \|^2_{-s}+ \|\nu_i (u) 
\|^2_{-s-1/2,\partial D}\right) ^{1/2}
$$
by $H ^{-s}_{\tau_i}(D,E_i)$ and $H ^{-s}_{\nu_i}(D,E_i)$ respectively.

\begin{corollary} \label{c.equiv.strong.taui}  
Let the differential complex $\{A_i\}$ be elliptic. Then the linear spaces 
$H^{-s}_{A_i}(D, E_i) $, $H ^{-s}_{\tilde \tau_i }(D,E_i)$  and $H 
^{-s}_{\tau_i }(D,E_i)$ coincide and their norms are equivalent. 
\end{corollary}

{\bf Proof.} The equivalence of the norms $\|\cdot \|_{-s,A_i}$ and  $\|\cdot 
\|_{-s,\tilde \tau_i}$ follows Theorem \ref{t.equiv.strong.tau}. Finally, as 
the complex $\{A_i\}$ is elliptic then Lemma \ref{l.projections} implies the 
equivalence of the norms $\|\cdot \|_{-s,\tilde \tau_i}$  and $\|\cdot 
\|_{-s,\tau_i}$. 
\hfill $\square$

\begin{corollary} \label{c.equiv.strong.nui}  
Let the complex $\{A_i\}$ be elliptic. Then linear spaces $H^{-s}_{A^*_{i-1}}(D,
E_i) $, $H ^{-s}_{\tilde \nu_i }(D,E_i)$  and $H ^{-s}_{\nu_i }(D,E_i)$ 
coincide and their norms are equivalent. 
\end{corollary}

{\bf Proof.} As the complex $\{A_i\}$ is elliptic then the complex 
$\{A^*_i\}$ is elliptic too. That is why Corollary \ref{c.equiv.strong.taui} 
implies the desired statement.  
\hfill $\square$

\begin{corollary} \label{c.equiv.strong.ti}  
If the complex $\{A_i\}$ is elliptic then the linear spaces $H^{-s}_{A_i \oplus 
A^*_{i-1}}(D, E_i) $ and  $H ^{-s}_{t}(D,E_i)$ coincide and their norms are 
equivalent. 
\end{corollary}

{\bf Proof.} As the complex $\{A_i\}$ is elliptic then the operator 
$A_i \oplus A^*_{i-1}$ has the injective principal symbol. Hence the statement 
follows from Theorem \ref{t.equiv.strong.tau}. 
\hfill $\square$

\begin{corollary} \label{c.identities}  
If the complex $\{A_i\}$ is elliptic then the following identities hold true:
$$
H^{-s}_{A_i \oplus A^*_{i-1}}(D, E_i) =
H^{-s}_{A_i}(D, E_i)\cap H^{-s}_{A^*_{i-1}}(D, E_i),
$$ 
$$
H ^{-s}_{t}(D,E_i) = H ^{-s}_{\tau_i }(D,E_i) \cap 
H ^{-s}_{\nu_i }(D,E_i). 
$$ 
\end{corollary}

\subsection{Weak boundary values of the tangential and normal parts}

Consider now the weak extension of an operator $A$ on the scale $H^{-s}(D,E)$. 
Namely, denote $H ^s_{A, w}(D,E)$ the set of the sections $u$ from $H ^{-s}(D,E)$ 
such that there is a section $f \in H ^{-s}(D, F)$ satisfying $Au = f$ in $H 
^{-s}(D, F, |\cdot |_{-s-1})$ (in particular, in the sense of distributions in 
$D$). As the operator $A$ is linear, this set is linear too. Clearly, 
\begin{equation} \label{eq.inclusion}
H ^{-s}_{A}(D, E) \subset H ^{-s}_{A, w}(D, E).
\end{equation}
It is natural to expect that these spaces coincide (cf. \cite{Fri}); we will 
prove it later. 

According to Corollary \ref{c.equiv.strong.taui}, we have $\tau_i(u) \in 
H^{-s-1/2}(\partial D, E_i)$ for all sections $u\in H^{-s}_{A}(D,E_i)$. Let us 
clarify the situation with the traces of the elements from $H^{-s}_{A_i,w} (D,E_i)$ 
for an operator $A_i$ from an elliptic complex. 

To this end, define pairing $(u,v)$ for $u \in H ^{-s}(D,E)$, 
$v \in C^{\infty} (\overline D,E)$ as follows. By the definition, one can find 
such a sequence $\{ u_{\nu} \}$ in $C^{\infty} (\overline D,E)$ that 
$\| u_{\nu} - u \|_{-s} \to 0$ if $\nu \to \infty$. Then 
$$
|(u_{\nu}-u_{\mu}, v)_{L^2 (D,E)}| \leq \| u_{\nu}-u_{\mu} \|_{-s} \| v 
\|_{H^{s}(D,E)}  \to  0 \mbox{ as }  \mu, \nu \to \infty. 
$$
Set $(u, v) = \lim_{\nu \to \infty}\limits (u_{\nu}, v)_{L^2 (D,E)}$.
It is clear that the limit does not depend on the choice of the sequence 
$\{ u_{\nu} \}$, for if $\| u_{\nu} \|_{-s} \to 0$, $\nu \to \infty$,
then $|(u_{\nu}, v)_{L^2 (D,E)}| \leq \| u_{\nu} \|_{-s} \| v \|_{H^{s}(D,E)}$
tends to zero too. This implies that for $u \in H ^{-s}(D ,E)$ and 
$v \in C^{\infty} (\overline D,E)$ we have the inequality:
$|(u, v)| \leq \| u \|_{-s} \| v \|_{H^{s}(D,E)}$.
Set $H (D, E) =\cup _{s=0}^\infty  H ^{-s} (D,E)$. Easily, 
the pairing $(u,v)_D$ is correctly defined for $u \in H (D, E)$ and
$v \in C^{\infty}(\overline D, E)$. The unions $\cup _{s=1}^\infty H ^{-s} 
(D,E)$ and $\cup _{s=1}^\infty H^{-s}_{A,w} (D, E)$ we denote by $H(D,E)$ and 
$H_{A} (D,E)$ respectively. 

As before, let $\Gamma$ be an open (in the topology of $\partial D$) connected 
subset of $\partial D$. The following definition is induced by 
(\ref{eq.Green.i}).
\begin{definition} \label{def.BV.tau} 
Let alone the correctness of this definition, we say that a 
distribution-section $u \in H_{A_i} (D, E_i)$, satisfying $A_i u =f$ in $D$ 
with $f \in H (D, E_{i+1})$,  has a weak boundary value $\tau_{i,\Gamma} ^w 
(u) = \tau_i (u_0)$ on $\Gamma$ for $u_0\in  {\mathcal D}'(\Gamma,E_i)$  if 
$$ 
(f,g)_D - (u,A_i^*g)_D =  \langle \star \tilde \nu_{i+1} (g), \tau_i(u_0)
\rangle_{\Gamma} \mbox{ for all }  g \in   C^\infty _{comp}(D\cup \Gamma, 
E_{i+1}) .
$$
\end{definition}

Formulae (\ref{eq.Green1}), (\ref{eq.Green.op}) and Theorem \ref
{t.equiv.strong.tau} imply that any section $u \in H ^{-s}_{A_i}(D,E_i)$ has 
a weak boundary value of the tangential part $\tau^w_{i,\partial D}(u)$ on 
$\partial D$ coinciding with the trace $\tau_{i,-s}(u) \in H^{-s-1/2} (\partial 
D, E_i)$. We are to connect the weak boundary values of the tangential parts 
with the so-called  {\it limit boundary values} of the solutions of finite orders 
of growth near $\partial D$ to elliptic systems (see \cite{Stra}, \cite{Roj}, 
\cite{Tark37}). Recall that a solution $u \in S_A (D)$ of an elliptic system 
$A$ has a finite order of growth near $\partial D$ if for any point  $x^0 \in 
\partial D$ there are a ball $B(x^0 ,R)$ and constants $c>0$, $\gamma >0$ such 
that 
$$
|v (x)| \leq c \,\, dist(x,\partial D)^{-\gamma} \mbox{ for all }
x \in B(x^0,R)\cap D.
$$
As $\partial D$ is compact, the constants $c$ and $\gamma$ may be chosen in 
such a way that this estimate is valid for all $x^0\in \partial D$. The space 
of all the solutions to $A$ of finite order of growth near $\partial D$ 
will be denoted $S^F_A (D)$.

Further, set $D_\varepsilon = \left\{x\in D: \rho(x)<-\varepsilon\right\}$. 
Then, for sufficiently small $\varepsilon>0$, the sets $D_\varepsilon\Subset D 
\Subset D_{-\varepsilon}$ are domains with smooth boundaries  $\partial 
D_{\pm\varepsilon}$ of class $C^\infty$. Besides, the vectors $\mp \varepsilon \nu (x)$ 
belong to  $\partial D_{\pm \varepsilon}$ for every $x \in \partial D$ (here 
$\nu (x)$ is the external normal unit  vector to the hyper-surface $\partial D$ 
at the point $x$). According to \cite[Theorem 9.4.7]{Tark37}, \cite{Roj}, if 
$A$ is elliptic and it satisfies the Uniqueness Condition \ref{US} then any 
solution $w \in S^F_{A^*A}(D)$ had the weak limit value $w^0 \in (C^\infty_{comp}
(\Gamma,E))^\prime$ on  $\Gamma$, i.e. 
$$ 
<w^0 , v >  = \lim_{\varepsilon \to +0} \int_{\partial D} 
 v (y) w (y-\varepsilon \nu (y))  ds (y) 
\mbox{ for all }  v \in C^\infty_{comp} (\Gamma, E).
$$

\begin{theorem} \label{t.BV} Let ${A_i}$ be an elliptic complex such that 
the operators $A_i \oplus A^*_{i-1}$, $0\leq i \leq N$, satisfy the Uniqueness 
Condition \ref{US}. Then every section $u \in H^{-s} _{A_i,w}(D,E_i)$  has 
the weak boundary value $\tau^w_{i,\partial D}(u)\in H^{-s-1/2}(\partial D,E_{i})$ 
in the sense of Definition \ref{def.BV.tau}, coinciding with the limit boundary 
value $\tau _i (w)$ of the solution $w=(u - {\mathcal G}_{\Delta_i} ^{(D,-s)}
A_{i}^* f - A_{i-1}{\mathcal G} ^{(D,-s+1)}_{\Delta_{i-1}} A_{i-1}^*u )$ from 
$S_{\Delta_i} ^F (D) $;  besides, $\tau^w_{i,\partial D}(u)$ does not depend 
on the choice of $f\in H^{-s-1}(D,E_{i+1}) $ with $A_i u =f $ in $D$. 
\end{theorem} 

{\bf Proof.} First of all we note that Lemma \ref{l.A}, Theorem \ref
{t.equiv.strong.tau} and Remark \ref{r.Dirichlet}, imply that  the operator 
${\mathcal G}^{(D,p)}_{\Delta_i} A^*_i$ continuously maps $H^{p-1}(D,E_{i+1})$ 
to $H ^p_{A_i}(D, E_i)$. Hence the sections  $w_1 = {\mathcal G}^{(D,-s)}
_{\Delta_i} A_{i}^*f \in H^{-s} _{A_i} (D,E_i)$ and $w_2={\mathcal G}^{(D,-s+1)}
_{\Delta_{i-1}} A ^* _{i-1} u \in H ^{-s+1}_{A_{i-1}}(D,E_{i-1})$ have the zero 
traces $t_{-s} (w_1)$ and $t_{-s+1} (w_2)$ on $\partial D$. In particular, 
$\tau_{i,-s} (w_1)=0$, $\tau_{i-1,-s+1} (w_2)=0$, and therefore $\tau ^w 
_{i,\partial D} (w_1)=0$,  $\tau ^w _{i-1,\partial D} (w_2)=0$. Besides, as 
$A_i \circ A_{i-1} \equiv 0$, we see that  $A_i (A_{i-1} w_2) =0$ in $D$ and 
$A_{i-1} w_2 \in H ^{-s}_{A_i, w}(D,E_i)$. According to Definition \ref
{def.BV.tau}, applied to $w_2$, we have: 
$$ 
(0,\psi)_D - (A_{i-1} w_2 ,A_{i}^* v)_D = - \langle \star \tilde \nu_{i} 
(A_i^* v),  \tau _{i-1}(w_2) \rangle _{\Gamma} + (  w_2,  A_{i-1}^* 
A_{i}^*   v)_D =0
$$
for all $v \in   C^\infty _{comp}(D\cup \Gamma, E_{i+1})$. Therefore $\tau 
^w _{i,\partial D} (A_{i-1} w_2)=0$ too.

It is clear now that the section $u\in H^{-s} _{A_i,w}(D, E_i)$ has the weak 
boundary value of $\tau ^w_{i,\partial D}(u)$ in the sense of Definition 
\ref{def.BV.tau} if and only if the section $w=(u - {\mathcal G} 
^{(D,-s)}_{\Delta_i} A_i ^*f -A_{i-1}  {\mathcal G} ^{(D,-s+1)}_{\Delta_{i-1}}  
A_{i-1} ^*u )$ has. By the construction $w \in H ^{-s}_{A_i ,w}(D,E_i)$ 
satisfies 
$$ 
\Delta_i  w = (A_{i} ^* A_{i} + A_{i-1} A_{i-1} ^*  ) u -  A_{i}  ^*f - 
A_{i-1} (A_{i-1} ^* u)  = 0 \mbox{ in } D. 
$$ 
In particular, this section belongs to $C^\infty (D,E_i)$, it has a finite 
order of growth near $\partial D$ (see \cite[Theorem 2.32]{ShT5}), and hence 
it has the limit boundary value $w^0\in(C^\infty_{comp}(\partial D,E_i))^\prime$  
on $\partial D$ (see \cite[Theorem 9.4.8]{Tark37}). Of course, the section  
$\tau _i (w^0) \in (C^\infty_{comp} (\partial D,E_i))^\prime$ is also defined 
because the function $\rho$ is of class $C^\infty$. Clearly, $\tau _i (w) =  
\tau_i (w^0)$ in the sense of the limit boundary values on $\partial D$.

As we have already noted, $w \in H ^{-s}_{A_i,w}(D,E_i)$ and 
$A_iw = f - A_i {\mathcal G}^{(D,-s)}_{\Delta_i} A^*_i f$  in $D$ where  
$(f - A_i  {\mathcal G}^{(D,-s)}_{\Delta_i} A^* _i f) \in H^{-s-1} (D,E_{i+1})$.
In particular, this means that 
$$ 
\langle \chi_D w, v \rangle = (w,v)_D \mbox{ for all } v\in C^{\infty} 
(\i{X},E_i),
$$
$$
\langle \chi _D (f - A_i {\mathcal G}^{(D,-s)}_{\Delta_i} A^*_i f), g \rangle 
= (f - A_i {\mathcal G}^{(D,-s)}_{\Delta_i} A^*_i f ,g )_D \mbox{ for all } g 
\in C^{\infty} (\i{X},E_{i+1}).
$$

Since the both $w$ and $A_iw$ are solutions to elliptic operators, i.e. 
$\Delta_i w= 0$  in  $D$, $\Delta_{i+1} (A_i w) =0$ in $D$ and 
they both have finite orders of growth near $\partial D$, then it follows 
from \cite[the proof of Theorem 9.4.7]{Tark37} that there is a sequence  of 
positive  numbers  $\{ \varepsilon_\nu\}$, tending to zero and such that 
$$ 
\langle \chi_Dw, v \rangle = \lim_{\varepsilon _\nu \to +0} 
\int_{D_{\varepsilon _\nu}} (w,v)_x dx 
\mbox{ for all } v\in C^{\infty}  (\i{X},E_i),
$$
$$ 
\langle \chi _D (f - A _i {\mathcal G}^{(D,-s)}_{\Delta_i} A^* _i f), g 
\rangle  = \lim_{\varepsilon _\nu \to +0} \int_{D_{\varepsilon _\nu}} (A_i 
w,g)_x dx \mbox{ for all } g\in C^{\infty}  (\i{X},E_{i+1}).
$$ 
By Whitney's Theorem, every smooth section over $\overline D$ may be extended 
up to a smooth section over $X$. Therefore 
$$ 
(w, v)_D  = \lim_{\varepsilon _\nu \to +0} \int_{D_{\varepsilon_\nu }} 
(w,v)_x dx \mbox{ for all } v\in C^\infty  (\overline D,E_i),
$$
$$ 
(f - A _i {\mathcal G}^{(D,-s)}_{\Delta_i} A_i^* f, g)_D = \lim_{\varepsilon 
_\nu\to +0} \int_{D_{\varepsilon_\nu}} (A_iw,g)_x dx \mbox{ for all } g\in 
C^\infty (\overline D,E_{i+1}).
$$ 
As $\tau_{i}({\mathcal G}^{(D,-s)}_{\Delta_i} A^*_i f + A_{i-1} {\mathcal 
G}^{(D,-s+1)}_{\Delta_{i-1}} A^*_i u) = 0$ on $\partial D$ in the sense of 
Definition \ref{def.BV.tau}, we see that Lemma \ref{l.projections}, formulae 
(\ref{eq.Green1}) and (\ref{eq.Green.i}) imply for all $g \in C^\infty  
(\overline D,E_{i+1})$:
$$  
(f, g)_D - (u, A^*_i g)_D = (f - A_i {\mathcal G}^{(D,-s)}_{\Delta_i} A^*_i f, 
g)_D - (w, A^*_ig)_D =
$$
$$ 
\lim_{\varepsilon _\nu \to +0} \left( 
\int_{D_{\varepsilon_\nu}} ((A_iw,g)_x -  (w,A_i^*g)_x) dx
\right) = 
$$ 
$$ \lim_{\varepsilon _\nu \to +0} \int_{\partial D_{\varepsilon_\nu}} 
 (\tau _i (w), \tilde \nu _{i+1} (g))_x \  ds (x) =  
\langle  \star  \tilde \nu_{i+1} (g), \tau_i ( w^0) \rangle _{\partial D} ,
$$
i.e. $\tau^w_{i,\partial D}(u) = \tau_i (w)$ on $\partial D$.
Now, if $\tilde f \in H ^{-s-1}(D,E_{i-1})$ satisfies $A_i u =\tilde f$ in 
$D$ then $\tilde w = (u - {\mathcal G}_{\Delta_i}^{(D,-s)} A^*_i \tilde f - 
A_{i-1}{\mathcal G}^{(D,-s+1)}_{\Delta_{i-1}} A^*_{i-1}u )$ and we have: $(w - 
\tilde w) = {\mathcal G}^{(D,-s)}_{\Delta_i} A^*_i (f - \tilde f) \in 
H^{-s}_{A_i} (D,E_i)$ with  $\tau ^w_{i,\partial D}(w - \tilde w) = 0$  on 
$\partial D$, i.e. the weak boundary value $\tau^w_{i,\partial D}(u)$ does not 
depend on the choice of the section $f \in H ^{-s-1}(D,E_{i+1})$ satisfying  
$A_iu =f$ in $D$.

Finally, we are to prove that the weak boundary value belongs to the 
corresponding Sobolev space $H ^{-s-1/2}(\partial D,E_{i})$. With this aim, 
fix a section $v_0 \in C^\infty (\partial D, E_{i+1})$. Then the section 
$g = {\mathcal P}_{\nabla_{E_{i+1}} ^*\nabla_{E_{i+1}}}^{(D)} \tilde \tau _i 
(v_0)$ (see the proof of Theorem \ref{t.equiv.strong.tau}) belongs to 
$C^\infty (\overline D, E_{i+1})$ and coincides with $\tilde \tau _i (v_0)$ on 
$\partial D$. Moreover, according to Remark \ref{r.Dirichlet} we have: 
\begin{equation} \label{eq.est}
\| g \|_{H^{s+1} (D,E_{i+1})} \leq \gamma_1 
\| \tilde \tau _i (v_0) \|_{H^{s+1/2} (\partial D, E_{i+1})} \leq \gamma_2 
\| v_0 \|_{H^{s+1/2} (\partial D, E_{i})}
\end{equation}
 with a positive constants $\gamma_1$, $\gamma_2$,  which does not depend on 
$g$ and $v_0$. Hence, by Definition \ref{def.BV.tau} and Lemma \ref
{l.projections}, we obtain: 
$$ 
|(\tau^w_{i,\partial D} (u), v_0)_{\partial D} |  = |\langle  \star \tilde 
\nu_{i+1} ( \tilde \tau_i (v_0)) , \tau^w_{i,\partial D} (u) \rangle_{\partial 
D} |= |\langle  \star \tilde \nu_{i+1} (g) , \tau^w_{i, \partial D} (u) \rangle
_{\partial D} |  =  
$$
$$
\left| (f,g)_D - (u,A_i^* g)_D\right| \leq \|f\|_{-s-1}  \|  g \|_{H^{s+1} 
(D,E_{i+1})} + \| u\|_{-s}  \| A^*_i g \|_{H^{s} (D,E_i)}.
$$  
As the map $A_i^*: H^{s+1} (D,E_{i+1}) \to H^{s} (D,E_i)$
is bounded, then the estimate implies that  
(\ref{eq.est}) 
$$ 
|(\tau^w_{i,\partial D} (u), v_0) | \leq \tilde \gamma (\| u\|_{-s} + 
\| f\|_{-s-1}) \| v _0 \|_{H^{s+1/2} (\partial D, E_{i+1})}  
$$
with a positive constant $\tilde \gamma $ which does not depend on 
$v_0$ and $u_0$. 

Hence, 
$$
\| \tau^w_{i,\partial D} (u) \|_{ H^{-s-1/2} (\partial D,E_{i} )} = 
\sup_{v \in C^{\infty} _{comp} ({\partial  D},E_i)}
\frac{|(\tau^w_{i,\partial D} (u) ,v ) _{\partial D} |}{\| v 
\|_{H^{s+1/2}({\partial D},E_i)}} \leq \tilde \gamma (\| u\|_{-s} +\| 
f\|_{-s-1}) . 
$$
Thus, the section $\tau^w_{i,\partial D}(u)$ belongs to the space 
$ H^{-s-1/2}(\partial D, E_i)$, which was to be proved. 
\hfill $\square$

\begin{corollary} \label{c.strong-weak.equiv} 
The spaces  $H^{-s}_{A_i} (D, E_i)$ and $H^{-s}_{A_i,w} (D,E_i)$ coincide.
\end{corollary} 

{\bf Proof.} Since (\ref{eq.inclusion}), it is enough to 
prove that $H^{-s}_{A_i,w}(D, E_i) \subset H ^{-s}_{A_i}(D, E_i)$. Fix a 
section $u \in H^{-s}_{A_i ,w}(D, E_i)$. Proving Theorem \ref{t.BV} we have 
seen that there is  $w \in S_{\Delta_i} ^F (D) \cap H^{-s} (D, E_i)$, 
satisfying 
$$ 
u = w + {\mathcal G}^{(D,-s)}_{\Delta_i} A_i^* f +
A_{i-1} {\mathcal G}^{(D,-s+1)}_{\Delta_{i-1}} A_{i-1}^* u.
$$
According to Remark \ref{r.Dirichlet}, the section $w$ is presented via its 
boundary values on $\partial D$ by the Poisson type integral $w= {\mathcal 
P}_{\Delta_i} ^{(D,-s)} t_i(w)$. Hence $w \in H^{-s}_t (D,E_i)$. Besides, 
Remark \ref{r.Dirichlet} imply that $w_1 = {\mathcal G}^{(D,-s)}_{\Delta_i} 
A_{i}^* f$  belongs to $H^{-s}_t (D, E_i)$ too. Thus, it follows from 
Corollary \ref{c.equiv.strong.ti} that the sections $w$ and $w_1$ belong 
$H^{-s}_{A_i \oplus A_{i-1}^*} (D, E_i) \subset H^{-s}_{A_i} (D, E_i)$. 

Take a sequence $\{ u_\nu\} \subset C^\infty (\overline D, E_i)$ approximating 
$u$ in the space $H^{-s} (D, E_i)$. It follows from  Remark \ref{r.Dirichlet} 
and \ref{l.A} that the sequence $\{ A_{i-1} {\mathcal G}^{(D,-s+1)}
_{\Delta_{i-1}} A_{i-1}^* u_\nu \} \subset C^\infty (\overline D, E_i)$ 
converges to  $A _{i-1} {\mathcal G}^{(D,-s+1)}_{\Delta_{i-1}} A^*_{i-1} u$ in 
the space $H^{-s}  (D,E_i)$. Moreover, $\{ A_{i}(A_{i-1} {\mathcal G}
^{(D,-s+1)}_{\Delta_{i-1}} A^*_{i-1} u_\nu) \equiv 0\} \subset C^\infty 
(\overline D, E_i)$ converges to zero in the space $H^{-s-1} (D,E_{i+1})$. 
Therefore  $A_{i-1} {\mathcal G}^{(D,-s+1)}_{\Delta_{i-1}} A_{i-1}^* u$ 
belongs to $H^{-s}_{A_i} (D, E_i)$. That is why the section $u$ belongs to 
this space too.
\hfill $\square$

\begin{corollary} \label{c.spaces} 
The differential operator $A_i$ continuously maps $H^{-s}_{A_i} (D, E_i)$ to  
$H^{-s-1}_{A_{i+1}} (D,E_{i+1})$.
\end{corollary} 

Similarly defining the spaces $H ^{-s}_{A_{i-1} ^*, w}(D,E_i)$ and  $H^{-s}
_{A_i \oplus A_{i-1}^*, w}(D,E_i)$ we easily obtain the following statements.

\begin{corollary} \label{c.strong-weak.equiv.nu} 
The spaces $H^{-s}_{A_{i-1}^*} (D,E_i)$ and $H^{-s}_{A_{i-1}^*,w} (D,E_i)$ 
coincide.
\end{corollary} 

\begin{corollary} \label{c.strong-weak.equiv.t} 
The spaces $H^{-s}_{A_i \oplus A_{i-1}^*} (D,E_i)$ and $H^{-s}_{A_i \oplus 
A_{i-1}^*,w} (D,E_i)$ coincide.
\end{corollary} 

As we have seen above, the scale $\{ H ^{-s}_{A_i}(D, E_i)\}$ is suitable 
for stating the Cauchy problem for the elliptic first order complex $\{A_i\}$. 
In order to do this we 
need to choose a proper spaces for the boundary Cauchy data on a surface 
$\Gamma \subset \partial D$. As we are interesting in the case $\Gamma \ne 
\partial D$, we will use one more type of the Sobolev spaces: the Sobolev spaces on 
closed sets (see, for instance, \cite[\S 1.1.3]{Tark37}). Namely,
let $H^{-s-1/2} (\overline \Gamma,E_i)$ stand for the factor space of  
$H^{-s-1/2} (\partial D,E_i)$ over the subspace of functions vanishing on a 
neighborhood of $\overline \Gamma$. Of course, it is not so easy to handle 
this space, but its every element extends from $\Gamma$ up to an element of 
$H^{-s-1/2} (\partial D,E_i)$. Further characteristic of this space may be 
found in \cite[Lemma 12.3.2]{Tark37}). We only note that if $\Gamma$ has 
$C^\infty$-smooth boundary (on $\partial D$), then 
$$
H^{-s-1/2}(\Gamma, E_i) \hookrightarrow  H^{-s-1/2} (\overline \Gamma, E_i)
\hookrightarrow  \tilde H^{-s-1/2} (\Gamma,E_i) .
$$ 
\begin{corollary} \label{c.BV.w} For every section  $u \in H ^{-s}_{A_i}(D, 
E_i)$ and every $\Gamma \subset \partial D$ there is the boundary value  
$\tau_{i,\Gamma}(u)$ in the sense of Definition  \ref{def.BV.tau}, belonging 
to $ H^{-s-1/2}(\overline \Gamma, E_i)$.
\end{corollary} 

As $\partial D$ is compact, $\cup _{s=1}^\infty H ^{-s-1/2} (\partial D, E_i) = 
{\mathcal D}^\prime (\partial D,E_i)$.  Set $\cup _{s=1}^\infty H^{-s-1/2} 
(\overline \Gamma,E_i) = {\mathcal D}^\prime (\overline \Gamma, E_i)$. Now 
Corollary \ref{c.strong-weak.equiv} immediately implies the following 
statements.

\begin{corollary} \label{c.BV.ww} For every $u \in H _{A_i}(D,E_i)$ and every 
$\Gamma \subset \partial D$ there is the boundary value $\tau_{i,\Gamma }(u)$ 
in the sense of Definition \ref{def.BV.tau}, 
belonging to $ {\mathcal D}^\prime (\overline \Gamma, E_i)$.
\end{corollary} 

\section{A homotopy formula}

In this section we will obtain an integral formula for  elements of 
the Soblev spaces with non-negative smoothness. Of course, for 
sufficiently smooth sections such formulae are well known (see, for instance, 
\cite[\S 2.4]{Ta1}).

From now on we additionally assume that the operators $\Delta_i$, $0\leq i\leq 
N$,  satisfy the Uniqueness Condition \ref{US}. Then each of these operators 
has a bilateral pseudo-differential fundamental solution, say, $\Phi _i$, 
on $\i{X}$ (see, for example, \cite[\S 4.4.2]{Tark37}). Schwartz kernel  
of the operator $\Phi_i$ is denoted by $\Phi _i (x,y)$, $x\ne y$. It is known, 
that $\Phi _i (x,y) \in C^\infty ((E_i \otimes E_i^*)\setminus \{x=y\})$ (see, 
for instance, \cite[\S 5]{Ta1}).  

For a section $f \in C^\infty (\overline D, E_{i+1}) $
we denote by $T_{i} f$ the following volume potential:
$$ T_{i}f (x)   =  (\Phi_i A_{i}^* \chi _D f) (x) = 
\int_D \langle (A_{i}^*)_y^T \Phi _i (x, \cdot) , f \rangle _y  dy . $$
 If $\partial D$ is smooth enough (e.g. $\partial D \in C^\infty$) 
then the potential $T_{i}$ induces a bounded linear operator
$$T_i: H^{s-1} (D, E_{i+1}) \to  H^{s} (D,E_{i}), \qquad s \in {\mathbb N}$$
(see, for example, \cite[1.2.3.5]{ReSch}). 

\begin{lemma} \label{l.T} For any domain $\Omega \Subset \i{X}$ 
with $\partial \Omega \in C^\infty$ the potential $T_{i}$ induces 
a bounded linear operator
$$
T_{i, \Omega}: H^{-s} (D, E_{i+1}) \to  H^{-s+1} _{A_i}
(\Omega, E_i), \quad s \in {\mathbb N}.
$$ 
Moreover for every section $f \in H^{-s}(D, E_{i+1})$ it is true that
$\Delta_i T_{i,\Omega} f = A_i^* \chi_D f$ in $\Omega \setminus \overline D$.  
\end{lemma}

{\bf Proof.} First of all we note that any smoothing operator 
$\tilde K$  of type $E_{i+1}\to E_i$ on $\i{X}$ induces for any $p$ a bounded 
linear operator 
$$
\tilde K \chi_D : H ^{-s} (D, E_{i+1}) \to  C ^p (\overline \Omega,E_i).
$$
As any two fundamental solutions differ on a smoothing operator, we may assume 
that $\Phi_i = {\mathcal G}_{\Delta_i}^{(X)}$. The principal advantage of 
${\mathcal G}_{\Delta_i}^{(X)}$ is in the fact that  the volume potential is 
$L^2(X,E_i)$-self-adjoint (see, for instance, \cite[formula (2.75)]{ShT5}). Besides, it 
has the transmission 
property (see \cite[\S 2.2.2]{ReSch}) and hence it continuously acts 
on the Sobolev scale:
$$
{\mathcal G}_{\Delta_i}^{(X)} \chi_D : H^{s-1} (D,E_i) \to H^{s+1} 
(\Omega,E_i), \quad {\mathcal G}_{\Delta_i}^{(X)} A_i^* \chi_D  : H^{s-1} 
(D,E_{i+1}) \to H^{s} (\Omega,E_i), \quad s \in {\mathbb N}.
$$
In particular, ${\mathcal G}_{\Delta_i}^{(X)} \chi_{\Omega} v$ belongs to $H^2 
_{loc}(\i{X}, E_i) \cap C^\infty (\overline \Omega, E_i)$ for all $v\in 
C^{\infty}(\overline \Omega , E_i)$ and, similarly,  ${\mathcal G}
_{\Delta_i}^{(X)} A_i^* \chi_{\Omega} g$ belongs to $H^1_{loc} (\i{X}, E_i) 
\cap C^\infty (\overline \Omega, E_i)$ for all $g \in C^{\infty}(\overline 
\Omega, E_{i+1})$. 
Then for all  $f \in C^\infty(\overline D, 
E_{i+1})$, $v \in C^{\infty} (\overline \Omega , E_i)$, $g \in C^{\infty} 
(\overline \Omega , E_{i+})$ we have: 
$$ 
(T_{i} f , v) _{\Omega}= ( {\mathcal G}_{\Delta_i}^{(X)}  A_i ^* \chi_{D}f,
\chi_{\Omega} v)_{X} = ( \chi_{D} f , A_i {\mathcal G}_{\Delta_i}^{(X)}  
\chi_{\Omega} v)_{X},
$$
$$   ( A_i T_{i} f , g) _{\Omega}= ( A_i  {\mathcal G}_{\Delta_i}^{(X)}  A_i ^* 
\chi_{D} f , \chi_{\Omega} g)_{X} = ( \chi_{D} f , A_i  {\mathcal 
G}_{\Delta_i}^{(X)}A_i ^*  \chi_{\Omega} g)_{X}.
$$
Therefore, we have 
\begin{equation} \label{eq.est.T}
 \|T _{i} f \|_{-s, A_i  ,\Omega} \leq C_1\, \|f \|_{-s-1,D} \mbox{ 
for all } f \in C^\infty(\overline D, E_{i+1} ), 
 \end{equation}
 \begin{equation} \label{eq.est.AT}
 \| A_i T _{i} f \|_{-s-1, A_i  ,\Omega} \leq C_2\, \|f \|_{-s-1,D} \mbox{ 
for all } f \in C^\infty(\overline D, E_{i+1} ), 
 \end{equation}
with positive constants $C_1$, $C_2$ do not depending on $f$.

Let now $f \in H^{-s-1}(D, E^{i+1})$. Then there is a sequence  
$  \{f_\nu \} \subset C^\infty(\overline D, E_{i+1})$ converging to 
$f$ in $H^{-s-1} (D, E_{i+1})$. According to (\ref{eq.est.T}), (\ref
{eq.est.AT}) the sequence  $\{ T _{i} f_\nu \}$ is fundamental in 
the space $H ^{-s}_{A_i}(\Omega, E_i)$; its limit we denote $T _{i,\Omega} f$. 
It is easy to understand that this limit does not depend on the choice 
of the sequence $\{f_\nu \}$ converging to $f$, and the estimates (\ref
{eq.est.T}), (\ref{eq.est.AT}) guarantee  that the operator $T_{i,\Omega}$, 
defined in this way, is bounded. Moreover, the properties of the fundamental 
solutions $\Phi_i$ means that each of the potentials 
$T _{i} f_\nu$ satisfies 
$$
(T_{i} f _\nu , \Delta_i v )_{\Omega} = 2 \langle A_i^* \chi_D f_\nu, 
v\rangle  =  ( \chi_D f_\nu, A_i v)_\Omega  \mbox{ for all }  
v \in C^\infty_{comp} (\Omega\setminus \overline D, E_i). 
$$ 
Passing to the limit with respect to $\nu \to \infty$ in the last equality we 
obtain the desired statement because the operators $\chi_D$ and $T_{i,\Omega}$ 
are continuous. 
\hfill $\square$

Further, for a section $v \in C^\infty (\overline D, E_i) $
we denote by $K_{i} f$ the following volume potential:
$$ 
K_{i}v  =  (\Phi_{i} A_{i-1} - A_{i-1}\Phi_{i-1} ) A_{i-1}^* \chi _D v . 
$$
Again, by the definition, it is a zero order pseudo-differential operator 
with the transmission property. If $\partial D$ is smooth enough (e.g. 
$\partial D \in C^\infty$) then the potential $K_{i}$ induces a bounded linear 
operator 
$$
K_i: H^{s} (D, E_{i}) \to  H^{s} (D,E_{i}), \qquad s \in {\mathbb Z}_+
$$
(see, for example, \cite[1.2.3.5]{ReSch}). 

\begin{lemma} \label{l.K} For any domain $\Omega \Subset \i{X}$ 
with $\partial \Omega \in C^\infty$ the operator $K_{i}$ induces 
a smoothing operator on $\overline \Omega$. In particular, for all $s \in 
{\mathbb N}$, $p \in {\mathbb N}$, it is bounded linear operator
$$
K_{i, \Omega}: H^{-s} (D, E_{i}) \to  C^p (\overline \Omega, 
E_i)\cap S_{\Delta_i} (\Omega). 
$$ 
\end{lemma}

{\bf Proof.} Indeed, by the definition of the fundamental solution, 
$$ 
\Delta_i (\Phi_{i} A_{i-1} - A_{i-1}\Phi_{i-1} ) v= A_{i-1}v - A_{i-1} v
=0 \mbox{ for all } v\in C^\infty _{comp}(\i{X}, E_{i-1}).
$$
Therefore the pseudo-differential operator $(\Phi_{i} A_{i-1} - 
A_{i-1}\Phi_{i-1} )$ (of order $(-1)$ on $X$) is smoothing on compact subsets 
of $\i{X}$. Now the similar statements follows for $K_i$.
\hfill $\square$

For $x \not \in \partial D$ we denote  $M_i v_0(x)$ the following 
Green integral with a density $v_0\in C^{\infty} (\partial D,E_i)$:
\begin{equation} \label{eq.MB.def}
M _i v_0  (x) = -
\int_{\partial D} G_{A_i}(\star A_i  \star^{-1} \Phi _i (x,\cdot) , v_0) =
-\int_{\partial D} (\tau_i (v_0)  , \tilde \nu_{i+1}  (A_i  
\star^{-1} \Phi _i (x,\cdot))_ y  \,ds (y),
\, x \not \in 
\partial D ;
\end{equation}
the last identity easily  follows from (\ref{eq.Green.i}). Thus we define the 
Green transform with a density $v_0 \in {\mathcal D}'(\partial D,E_i)$ as the result 
of the action of the distribution $v_0$ on the "test-function"{}  $(- \tilde 
\nu_i  (A_i  \star^{-1} \Phi _i (x,\cdot)) \in C^\infty (\partial D,E_i)$: 
$$
M _i v_0  (x) = - (v_0  , \tilde \nu_{i+1}  (A_i  \star^{-1} \Phi _i (x,\cdot))
_{\partial D} =- (\tau_i (v_0)  , \tilde \nu_{i+1}  (A_i  \star^{-1} \Phi _i 
(x,\cdot))_{\partial D}, \qquad x \not \in \partial D.
$$
By the construction, $M _i v_0 \in S_{\Delta_i} (\i{X} \setminus 
\mbox{supp} \ v_0, E_i)$ as a parameter dependent distribution; here 
$\mbox{supp} \ v_0$ is the support of $v_0$.

Again, if $\partial D$ is smooth enough (e.g. $\partial D \in C^\infty$) then 
the potential $M_{i}$ induces a bounded linear operator 
$$
M_i: H^{s-1/2} (\partial D, E_{i}) \to  H^{s} (D,E_{i}), \qquad s \in 
{\mathbb N}
$$
(see, for example, \cite[1.2.3.5]{ReSch}). 

Now using Stokes formula and the potentials $T_i$, $M_i$, $K_i$ we arrive 
to a homotopy formula for the complex $\{A_i\}$ and sections $u \in C^\infty 
(\overline D,E_i)$ (see \cite[Theorem 2.4.8]{Ta1}):
\begin{equation} \label{eq.MB.bas}
M _{i} u  + T  _{i} A_i u + A_{i-1}
T  _{i-1}  u  + K_{i} u = \chi_D u.
\end{equation}
Of course, the continuity of the operators $T_i$, $M_i$, $K_i$  on the Sobolev
spaces implies that formula (\ref{eq.MB.bas}) is still valid for all the 
sections $u \in H^s (D,E_i)$, $s \in \mathbb N$. We are to extend the 
homotopy formula for the complex $\{A_i\}$ on the scale $H^{-s}_{A_i}(D, 
E_i)$, $s \in {\mathbb Z}_+$.

\begin{lemma} \label{l.MB}  For any domain $\Omega \Subset \i{X}$ such that  
$\partial \Omega \in C^\infty$ and $D\subset \Omega$ the potential 
$M$ induces bounded linear operators 
$$
M_{i,D}: H ^{-s-1/2} (\partial D, E_i) \to  H ^{-s}_{A_i}(D,E_i), \quad 
M_{i,\Omega}:  H ^{-s-1/2} (\partial D, E_i) \to  H ^{-s}  (\Omega, E_i). 
$$
\end{lemma}

{\bf Proof.} As we already have seen above (see Remark \ref{r.Dirichlet} and 
Corollary \ref{c.equiv.strong.ti}), for every section $ v ^0 \in H^{-s-1/2} 
(\partial D, E_i)$ the Poisson integral ${\mathcal P}_{\Delta_i}^{(D)} v ^0  
\in H ^{-s}_{A_i \oplus A_{i-1}^*}(D, E_i)$ satisfies $t_i( {\mathcal 
P}_{\Delta_i}^{(D)} v ^0) = v^0 $. Set 
$$ 
M _{i,D} = (I - T  _{i,D} A_i - A_{i-1}  T  _{i, D} - K_{i,D}) \ {\mathcal 
P}_{\Delta_i} ^{(D)}: H^{-s-1/2} (\partial D, E_i) \to H ^{-s} _{A_i}(D, E_i),
$$
$$ 
M _{i,\Omega} = (\chi_{D} - T  _{i,\Omega} A_i - A_{i-1}  T  _{i, \Omega} - 
K_{i,\Omega}) \ {\mathcal P}_{\Delta_i} ^{(D)}: H^{-s-1/2} (\partial D, E_i) 
\to H ^{-s} (\Omega, E_i).
$$
It follows from Lemmas \ref{l.T}, \ref{l.K} and the continuity of the 
operators ${\mathcal P}_{\Delta_i} ^{(D)}$ and $\chi_{D}$ that the defined above 
operators $M _{i,D}$, $M_{i,\Omega}$ are bounded. Let us see that $M_{i,D}$ and 
$M_{i,\Omega}$ coincide with $M_i$ on $C^{\infty} (\partial D, E_i)$. Indeed, 
if $v^0 \in C^{\infty} (\partial D, E_i)$ then Remark \ref{r.Dirichlet} 
implies ${\mathcal P}_{\Delta i} ^{(D)} v^0  \in C^{\infty} (\overline D, E_i)$
and  
$$
M_iv^0 = M_i {\mathcal P}_{\Delta_i} ^{(D)}v^0 = M _i\tau_i({\mathcal P}_
{\Delta_i} ^{(D)}v^0).
$$
Now using a homotopy formula (\ref{eq.MB.bas}) we obtain:
$$ 
\chi_{D} {\mathcal P}_{\Delta_i} ^{(D,i)}  v^0  = M _i v^0  + T_{i,D} A_i 
{\mathcal P}_{\Delta_i} ^{(D)} v^0 +  A_{i-1} T_{i-1, D}  {\mathcal 
P}_{\Delta_i} ^{(D)} v^0 + K_i  {\mathcal P}_{\Delta_i} ^{(D)} v^0 .
$$
Since $C^{\infty} (\partial D, E_i)$ is dense in 
$H^{-s-1/2} (\partial D, E_i)$ then $M_i$ continuously extends from 
$C^{\infty} (\partial D, E_i)$ onto $H^{-s-1/2} (\partial D, E_i)$  
as defined above operators $M_{i,D}$, $M_{i,\Omega}$. Moreover, it is easy 
to understand that the sections $M_{i,D} v^0 $, $M_{i,\Omega} v^0$ are coincide 
with the distributions $Mv^0$ on $D$ and $\Omega \setminus \mbox{supp} \, v^0 $
respectively. 
\hfill $\square$

\begin{theorem} \label{t.MB} 
For every section $u \in H_{A_i}(D, E_i)$ the following formulae hold:
\begin{equation} \label{eq.MB-}
M _{i,D} u  + T  _{i, D} A_i u + A_{i-1}
T  _{i-1, D}  u  + K_{i,D} u = u, 
\end{equation}
\begin{equation} \label{eq.MB}
M _{i, \Omega} u  + T _{i,\Omega } A_i  u  
 + A_{i-1} T  _{i-1, \Omega}  u + K_{i,\Omega} u   = \chi_{D} u. 
\end{equation} 
\end{theorem}  

{\bf Proof.} Pick $u\in H_{A_i}(D, E_i)$. Then $u  \in H^{-s}_{A_i}(D, E_i)$ 
with a number $s\in {\mathbb Z}_+$ and there is $\{u_\nu \}\subset C^\infty
(\overline D, E_i)$ converging to $u$ in the space $H^{-s}_{A_i} (D,E_i)$. Now 
the homotopy formula (\ref{eq.MB.bas}) implies
\begin{equation} \label{eq.MB.1}
M_i u_\nu  + T _{i} A_i u_\nu + A_{i-1} T _{i}  u_\nu + K_i u_\nu 
 = \chi_{D} u_\nu. 
\end{equation}
Passing to the limit in the spaces $H^{-s} 
_{A_i}(D,E_i)$ and  $H ^{-s}(\Omega, E_i)$  
with respect to $\nu \to \infty$ in (\ref{eq.MB.1}) we obtain 
(\ref{eq.MB-}) and (\ref{eq.MB}) respectively because of 
Lemmas \ref{l.T}, \ref{l.K}, \ref{l.MB}.
\hfill $\square$

\begin{remark} \label{r.T} 
{\rm 
Let $f \in H^{-s-1}(D, E_{i+1})$. If $\Omega$, $\Omega_1$ are bounded domains 
in $\i{X}$ (with smooth boundaries) containing $D$ then sections 
$T_{i,\Omega} f \in H^{-s}(\Omega , E_i)$ and $T_{i,\Omega_1} f \in 
H^{-s}(\Omega_1, E_i) $ belong to  $S_{\Delta_i} (\Omega\setminus \overline D)$
and $S_{\Delta_i} (\Omega_1\setminus \overline D)$ respectively. Since they 
are constructed as the limits of  the same sequence of sections converging in different 
spaces, they coincide in $(\Omega_1 \cap \Omega) \setminus \overline D$.
The same conclusion is obviously valid for the smoothing operators 
$K_{i,\Omega}$ and $K_{i,\Omega_1}$. Moreover,  as the operators $M_{i,\Omega}$ 
and $M_{i,\Omega_1}$ are constructed with the use of $T_{i,\Omega}$, 
$K_{i,\Omega}$ and $T_{i,\Omega_1}$, $K_{i,\Omega_1}$ respectively, 
this is also true for the sections of the type $M_{i,\Omega} (v^0)$ 
with $v^0 \in  H^{s-1/2} (\partial D, E_i)$. Since $\Omega \subset \i{X}$ is 
arbitrary, the Uniqueness Condition \ref{US} allows us to say about the sections 
$T_{i} f$ and $M_i v^0 $ from  $S_{\Delta_i}^F(\i{X} \setminus \overline D)$
such that $T_i f = T_{i, \Omega} f \in H^{-s}(\Omega, E_i)$, $M_i v^0 = 
M_{i,\Omega} v^0 \in H^{-s}(\Omega , E_i)$ for any domain $\Omega \supset D$. 
\hfill $\square$
}
\end{remark}

\section{The Cauchy problem in spaces of distributions} 

\begin{problem} \label{pr.Cau} Given 
$u_0 \in {\mathcal D}^\prime (\overline \Gamma, E_{i})$, $f \in H_{A_{i+1}}(D, 
E_{i+1})$ find a section $u \in H_{A_i} (D,E_i)$ such that
\begin{equation*} 
A_i u = f \mbox{ in } D, \qquad 
\tau_i (u) = \tau_i (u_0) \mbox{ on } \Gamma,
\end{equation*}
in the sense of Definition \ref{def.BV.tau}, i.e. 
\begin{equation} \label{eq.Cau}
(u,A_i^*g)_D =  (f,g)_D  -\langle \star \tilde \nu_{i+1} (g), \tau_i(u_0)\rangle
_{\Gamma} \mbox{ for all }  g \in   C^\infty _{comp}(D\cup \Gamma, E_{i+1}) .
\end{equation}
\end{problem}
If $i=0$ then $A_0$ has the injective principal symbol and the Cauchy problem 
has no more than one solution (see, for instance, \cite[Theorem 10.3.5]
{Tark37}). Clearly it may have infinitely many solutions if $i>0$. Usually 
the Uniqueness Theorem of the Cauchy problem for $i>0$ is valid in 
co-homologies under some convexity conditions on $\partial D \setminus 
\Gamma$ (cf. \cite[Corollary 3.2]{NaSchTa}). Instead of looking for a version 
of Uniqueness Theorem we will try to choose a canonic solution of the Cauchy 
problem (see \S \ref{s.L2} below for solutions in $H^0_{A_i} (D,E_i)$). 

We easily see that $f$ and $u^0$ should be  coherent. Namely, as $A_i^* 
A^*_{i+1} \equiv 0$, taking $g= A^*_{i+1}w$ with $w \in C^\infty 
_{comp}(D\cup\Gamma, E_{i+2})$ in (\ref{eq.Cau}) we conclude 
that for the solvability of problem \ref{pr.Cau} it is necessary that 
\begin{equation} \label{eq.CR}
(f,A^*_{i+1}w)_D  = \langle \star \tilde \nu_{i+1} (A^*_{i+1} w),\tau_i(u_0)
\rangle _{\Gamma} \mbox{ for all } w \in C^\infty _{comp}(D\cup\Gamma, E_{i+2}).
\end{equation}

Let us discuss this. First we note that, due to Corollary \ref{c.spaces} and 
to the properties of the complex, $A_{i+1} f=0$ in $D$ if the Cauchy problem is 
solvable. This corresponds to $w \in C^\infty _{comp}(D, E_{i+2})$ in 
(\ref{eq.CR}).

Besides, the operator $A_i$ induces the tangential operator $\{ A_{i,\tau}\}$ on 
$\partial D$ (see, for instance, \cite[\S 3.1.5]{Ta1}). More precisely, let 
$\hat u^0 \in {\mathcal D}'(\partial D, E_i)$. Pick a section $\hat u \in H_{A_i} 
(D,E_i)$ satisfying $\tau_i (\hat u) = \tau_i (\hat u^0)$ on $\partial D$ 
(there is at least one such a section, ${\mathcal P}^{(D)}_{\Delta_i} \tau_i
(\hat u^0)$). Then set $ A_{i,\tau} \hat u^0 = \tau _{i+1} (A_i \hat u)$. If 
we fix $g \in   C^\infty(\partial D, E_{i+1})$ then, by Remark \ref
{r.Dirichlet}, the section $w={\mathcal P}^{(D)}_{\Delta_{i+2}} \tilde 
\tau_{i+1} (g)$ belongs to the space $C^\infty (\overline D, E_{i+2})$. Now, 
easily, Definition \ref{def.BV.tau} and Lemma \ref{l.projections} imply that 
$$
\langle \star g, A_{i,\tau} \hat u^0 \rangle = \langle \star \tilde \nu_{i+2} 
(\tilde \tau _{i+1} (g)), \tau _{i+1} (A_i \hat u) \rangle = 
\langle \star \tilde \nu_{i+2} (w) , \tau _{i+1} (A_i \hat u) \rangle =
$$ 
\begin{equation} \label{eq.Atau}
(A_{i} \hat u , A_{i+1} ^* w)_D = \langle  \star \tilde \nu_{i+1} (A_{i+1} ^* 
w), \tau_i(\hat u) \rangle = \langle  \star \tilde \nu_{i+1} (A_{i+1} ^* w), 
\tau_i(\hat u^0) \rangle .
\end{equation}
In particular, this means that $A_{i,\tau} \hat u^0 $ does not depend on the 
choice of $\hat u \in H_{A_i} (D,E_i)$ with $\tau_i (\hat u) = \tau_i (\hat 
u^0)$ on $\partial D$.

\begin{lemma} \label{l.CR} 
For the Cauchy data $u_0$ and $f$, identity {\rm (\ref{eq.CR})} holds if and 
only if $A_{i+1} f=0$ in $D$ and $\tau _{i+1,\Gamma}(f) = A_{i,\tau} u^0 $ on 
$\Gamma$. 
\end{lemma}

{\bf Proof.} Indeed, as we have noted above, (\ref{eq.CR}) implies $A_{i+1} 
f=0$ in $D$. Then, similarly to (\ref{eq.Atau}), it follows from Definition 
\ref{def.BV.tau} that, with $w={\mathcal P}^{(D)}_{\Delta_{i+2}} \tilde 
\tau_{i+1} (g)$, 
$$
\langle \star g,  \tau _{i+1} (f) \rangle = \langle \star \tilde \nu_{i+2} 
(\tilde \tau _{i+1} (g)), \tau _{i+1} (f) \rangle = 
\langle \star \tilde \nu_{i+2} (w) , \tau _{i+1} (f) \rangle =
(f , A_{i+1} ^* w)_D 
$$
for all $g \in   C^\infty(\partial D, E_{i+1})$ if $A_{i+1} f=0$ in $D$.  
Therefore taking $\hat u^0= u^0$ on $\Gamma$ and $g \in C^\infty_{comp}
(\Gamma, E_{i+1})$ in (\ref{eq.Atau}) we conclude that $\tau _{i+1,\Gamma}(f) 
= A_{i,\tau} u^0 $ on $\Gamma$ too, if identity (\ref{eq.CR}) holds.

Back, if $A_{i+1} f=0$ in $D$ and $\tau _{i+1,\Gamma}(f) = A_{i,\tau} u^0 $ on 
$\Gamma$ then , again applying  Definition \ref{def.BV.tau} and calculating 
as in (\ref{eq.Atau}), we obtain for all $w \in C^\infty _{comp}(D\cup\Gamma, 
E_{i+2})$:
$$
(f,A^*_{i+1}w)_D =\langle \star \tilde \nu_{i+2} (w) , \tau _{i+1} (f) \rangle 
= \langle \star \tilde \nu_{i+2} (w) , A_{i,\tau} u^0  \rangle = 
\langle  \star \tilde \nu_{i+1} (A_{i+1} ^* w), \tau_i( u^0) \rangle ,
$$
which was to be proved. \hfill $\square$

It is important to note that Lemma \ref{l.CR} allows the point wise check of necessary 
solvability conditions for Problem \ref{pr.Cau}, at least if 
the Cauchy data $f$ and $u^0$ are smooth.  

Now choose a domain $D^+$ in such a way that the set $\Omega = D\cup \Gamma 
\cup D^+$ is a bounded domain with smooth boundary in $\i{X}$. It is 
convenient to denote $F^\pm $ the restrictions of a section $F$ onto $D^\pm$ 
(here $D^- = D$). 

Further, for $u^0 \in H ^{-s-1/2} (\overline \Gamma, E_i)$, choose a 
representative $\tilde u^0 \in H ^{-s-1/2} (\partial D, E_i)$.  We have 
seen above the potentials $M_i \tilde u_0$ and $T_{i}f$ satisfy $\Delta_i (M_i 
\tilde u_0)= 0$ and $\Delta_i (T_i f)  = 0$  everywhere outside 
$\overline D$ as parameter dependent distributions. Hence the section
$$ 
F_i  = M_{i, \Omega} \tau_i ( \tilde u_0) + T_{i,\Omega} f 
$$
belongs to $S_{\Delta_i} (D^+) \cap H(\Omega,E_i)$. The Green formula (\ref 
{eq.MB}) shows  that the potential $F_i$ contains a lot of information on 
solvability conditions of Problem \ref{pr.Cau}.

Denote $\chi_D (H(D, E_i))$ the image of the space $H(D, E_i)$ under the map 
$\chi_D : H(D, E_i) \to 
H(\Omega, E_i) $ (see map (\ref{eq.chi})). 

\begin{theorem} \label{t.AKy}
Let $\Delta_{i-1}$, $\Delta_i$, $\Delta_{i+1}$ satisfy the Uniqueness 
Condition {\rm \ref{US}}. Then the Cauchy Problem {\rm \ref{pr.Cau}} is 
solvable if and only if condition {\rm (\ref{eq.CR})} holds true and there is 
a section ${\mathcal F }_i\in H(\Omega, E_i)$ such that $A_i \Delta_i {\mathcal F }_i 
= 0$ in $\Omega$ and $(F_i - {\mathcal F}_i) \in\chi_D (H(D,E_i))$.
\end{theorem}

{\bf Proof}. Let Problem \ref{pr.Cau} be solvable and $u$ be its
solution. The necessity of condition (\ref{eq.CR}) is already proved. Set
\begin{equation} \label{eq.phi.1}
{\mathcal F}_{i,u} = M_{i,\Omega} \tau_i (\tilde u^0) +  T_{i, \Omega} f - \chi_
{D} u. 
\end{equation}
Lemmas \ref{l.T}, \ref{l.T}, \ref{l.MB} and Remark \ref{r.T} imply that 
${\mathcal F}_{i,u} \in H^{-s}(\Omega, E_i)$ with some $s\in {\mathbb Z}_+$.
Clearly $(F_i - {\mathcal F}_i)= \chi_D u \in  \chi_D (H(D,E_i)) $.
Then it follows from homotopy formula (\ref{eq.MB}) that:
\begin{equation} \label{eq.F0}
{\mathcal F} _{i,u} = 
M_{i,\Omega}(  \tau_i  (\tilde u^0) -  \tau_i (u) )) - A_{i-1} T_{i-1, \Omega} 
u - K_i u.
\end{equation}
Since $( \tau_i (\tilde u^0) - \tau_i (u) )=0$ on $\Gamma$ then $M _{i,\Omega}( 
\tau_i(\tilde u^0) -  \tau_i (u))$ belongs to $S_{\Delta_i} (\i{X}\setminus
\Gamma)$ as a parameter dependent distribution. That is why, using Lemma 
\ref{l.T}, we obtain: 
\begin{equation} \label{eq.DeltaF}
\Delta_i {\mathcal F }_{i,u}  = - \Delta _i A_{i-1} T_{i-1, \Omega} u = -  A_{i-1}
\Delta_i T_{i-1, \Omega} u =  -  A_{i-1} A_{i-1} ^* \chi_D u \mbox{ in }\Omega .
\end{equation}
In particular, $A_i \Delta_i {\mathcal F }_{i,u} = 0$ in $\Omega$.

Back, let there be sections ${\mathcal F}_{i} \in H(\Omega,E_i) $ and $u \in H(D, 
E_i)$ such that $A_i \Delta_i {\mathcal F }_i = 0$ in $\Omega$ and 
\begin{equation} \label{eq.sol}
\chi_D u = F_i - {\mathcal F}_i. 
\end{equation}
Let us show that the section $u$ is a solution to Problem \ref{pr.Cau}. With
this aim we consider the following functional $w ( \tilde  u^0) $ 
on the space $C^\infty (\overline D,E_{i+1})$:
$$ 
\langle w ( \tilde  u^0), v \rangle = (\tau_i ( \tilde u^0), \tilde \nu _{i+1}  
(v))_{\partial D} \mbox{ for all } v \in C^\infty (\overline D, E_{i+1}).
$$
As $\tilde  u^0 \in {\mathcal D}^\prime (\partial D, E_i)$ then $\tilde  u^0 
\in H^{-s-1/2} (\partial D, E_i)$ with some $s \in {\mathbb Z}_+$ 
and hence for all  $v \in C^\infty (\overline D, E_{i+1})$ we have: 
$$ 
|\langle w ( \tilde  u^0) , v \rangle | \leq \| \tau_i( \tilde  u^0)\| 
_{-s-1/2, \partial D} \| \tilde \nu_{i+1} (v)  \| _{s+1/2, \partial D} \leq 
C \  \| \tau_i( \tilde  u^0)\| _{-s-1/2, \partial D} \| v  \| _{s+1, D}
$$ 
with a constant $C>0$ which does not depend on $\tilde  u^0$ and $v$. 
Therefore $w ( \tilde  u^0) \in H^{-s-1}(D, E_{i+1})$ and its support 
belongs to $\partial D$. 

Clearly, $C^\infty _{comp} (D\cup \Gamma, E_{i+1}) \subset C^\infty 
_{comp}(\Omega, E_{i+1})$ and Whitney theorem implies that every section from  
$C^\infty _{comp} (D\cup \Gamma, E_{i+1})$ may be extended up to an element of 
the space $C^\infty _{comp}(\Omega, E_{i+1})$. Thus,  (\ref{eq.Cau}) is 
equivalent to the following identity:
\begin{equation} \label{eq.Cau.distr}
g = A_i \chi_D u - \chi_D f + \chi_D w ( \tilde  u^0) \equiv 0 
\mbox{ in } \Omega. 
\end{equation}
That is why $u$ is a solution to Problem \ref{pr.Cau} if and only if $u \in 
H_{A_i}(D, E_i)$ and the identity (\ref{eq.Cau.distr}) holds. By the very 
construction, $g$ belongs to ${\mathcal D}' (\Omega, E_{i+1})$ and its support 
lies in $\overline D$.

Then for all $v \in C^\infty _{comp}(\Omega, E_{i+1})$ we have 
$$
\langle g , \Delta_{i+1} v \rangle _{\Omega} = (\chi_D u, A_i ^*  \Delta_{i+1} 
v)_\Omega-(\chi_D f, \Delta_{i+1} v)_\Omega + (\chi_D w(\tilde u^0),\Delta 
_{i+1} v)_{\Omega}  
$$
$$
(F_i - {\mathcal F}_i, \Delta_i A_i^*v)_\Omega - (f, \Delta _{i+1} v)_D + 
(\tau_i(\tilde u^0), \tilde \nu _{i+1} (\Delta_{i+1} v) ) _{\partial D} = 
$$
\begin{equation} \label{eq.AKy1}
(F_i, \Delta_i A_i ^* v)_\Omega - (f, \Delta _{i+1} v)_D + 
(\tau_i (\tilde u^0), \tilde \nu_{i+1} (\Delta _{i+1} v) ) _{\partial D} ,
\end{equation}
because $A^*_i \Delta_{i+1} = \Delta_{i} A^*_i $ and $A_i \Delta _i{\mathcal F }_i 
= 0$ in $\Omega$.

Further, by Lemma \ref{l.T}, we see that for all $v \in C^\infty _{comp}
(\Omega, E_{i+1})$, 
\begin{equation} \label{eq.AKy2}
(T_{i,\Omega}f, \Delta _i A_i ^*v)_\Omega = (A_i ^* \chi_D f,  A_i ^* v)_\Omega
 = (f,  A_i A_i^* v)_D . 
\end{equation}
Set $\tilde u = {\mathcal P}_{\Delta_i}^{(D)} \tau _i (\tilde u^0) $. This 
section belongs to $H_{A_i \oplus A_{i-1}^*} (D, E_i)$ (see Remark \ref
{r.Dirichlet} and Corollary \ref{c.equiv.strong.ti}). By the definition, 
$\tau_i(\tilde u ) = \tau _i (\tilde u_0) $ on $\partial D$. Now Lemma \ref
{l.MB}, the properties of the fundamental solutions and Definition \ref
{def.BV.tau} imply that for all $v \in C^\infty _{comp}(\Omega, E_{i+1})$ we 
have:
$$
(M_{i,\Omega} \tau _i(\tilde u^0), \Delta_i A_i^*v)_\Omega = (\chi_D \tilde u  
- T_{i,\Omega} A_i  u  - A_{i-1} T_{i-1,\Omega} \tilde u - K_{i,\Omega} u, 
\Delta _i A_i ^*v)_\Omega =
$$
\begin{equation} \label{eq.AKy3}
(\tilde u,  A_i ^* A_i A_i ^*v )_D - (A_i \tilde u , A_i A_i^*v )_D= - (\tau_i 
(\tilde u^0), \tilde \nu_{i+1} ( A_i A_i ^* v) _{\partial D} .
\end{equation}

Therefore, using (\ref{eq.AKy1}), (\ref{eq.AKy2}), (\ref{eq.AKy3}) we conclude 
that  
$$
\langle g,\Delta_{i+1} v\rangle _{\Omega} = -(f, A_{i+1} ^*  A_{i+1} v )_D +
(\tau_i (\tilde u^0), \tilde \nu_{i+1} ( A_{i+1} ^*  A_{i+1} v) ) _{\partial D} 
= 0
$$
for all $v \in C^\infty _{comp}(\Omega, E_{i+1})$ because of condition 
(\ref{eq.CR}).

Thus, $\Delta_{i+1} g =0$ in $\Omega$ and $g=0$ in $D^+$. It follows from  
Uniqueness Condition \ref{US} that $g \equiv 0$ in $\Omega$, i.e. identity 
(\ref{eq.Cau}) holds. In particular this means that  $A_i u=f$ in $D$ and, 
by Corollary \ref{c.strong-weak.equiv}, we see that $u\in 
H_{A_i} (D,E_i)$, which was to be proved.
\hfill $\square$

\begin{corollary} \label{c.AKy} Let $f \in H^{-s-1} (D, E_{i+1})$, 
$u^0 \in H^{-s-1/2}(\overline \Gamma, E_i)$.  The Cauchy problem 
{\rm \ref{pr.Cau}} is solvable in the space $H^{-s}_{A_i} (D, E_i)$ if and 
only if condition {\rm (\ref{eq.CR})} is fulfilled and there is a section 
${\mathcal F }_i \in H^{-s} (\Omega, E_i) $ satisfying $A_i \Delta_i {\mathcal F }_i 
=0$ in $\Omega$ and such that $(F_i - {\mathcal F}_i) \in  \chi_D (H^{-s}(D, 
E_i))$. 
\end{corollary}

{\bf Proof.} Indeed, if Problem \ref{pr.Cau} is solvable in $H_{A_i} ^{-s} (D, 
E_i)$, then condition {\rm (\ref{eq.CR})} is fulfilled  and ${\mathcal 
F}_i=F_i- \chi_D u$ (see (\ref{eq.phi.1})). Hence, by Lemma \ref{l.MB}, the 
section ${\mathcal F }$ belongs to $H^{-s} (\Omega, E_i)$ and  $(F_i - 
{\mathcal F}_i) \in  \chi_D (H^{-s}(D, E_i ))$. 

Back, if condition {\rm (\ref{eq.CR})} is fulfilled, ${\mathcal F }_i\in H 
^{-s}(\Omega, E_i) $ satisfies $A_i \Delta_i {\mathcal F }_i =0$ in $\Omega$ 
and $(F_i - {\mathcal F}_i) \in  \chi_D (H^{-s}(D, E_i)) $ then Problem \ref
{pr.Cau} is solvable. Besides, one of its solutions $u$ is given by formula (
\ref{eq.sol}). In particular, $\chi_D u = (F_i-{\mathcal F}_i)$ belongs to 
$H^{-s} (\Omega, E_i)$. Pick $v \in C^{\infty} (\overline D, E_i)$. Then, by 
Whitney Theorem, there is a section $V \in C^{\infty} (\overline \Omega, E_i)$ 
with $\|V\|_{s,\Omega}=\|v\|_{s,D}$ and $v=V$ in $D$. By the definition, 
$$ 
|(u,v)_D| = |(\chi_Du, V)_{\Omega}|\leq \| \chi_D u\|_{-s,\Omega}  \|v\|_{s,D},
$$
i.e. $u \in H^{-s} (D, E_i)$. Finally, as $A_i u =f \in H ^{-s-1}
(D, E_{i+1} ) $,  then $u \in H^{-s} _{A_i}(D, E_i)$ according to Corollary  
\ref{c.strong-weak.equiv}. 
\hfill $\square$

If $i=0$ then the operator $A_0$ has injective principal symbol and Theorem 
\ref{t.AKy} has the following form (cf. \cite{Tark37}, \cite{SheShlL2} for  
the operators with real analytic coefficients and $f=0$).

\begin{corollary} \label{c.AKy.ell} 
Let $f \in H (D, E_1)$, $u^0 \in {\mathcal D}' (\overline \Gamma, E_0)$.  
The Cauchy Problem {\rm \ref{pr.Cau}} is solvable in the space $H_{A_0} (D, 
E_0 )$ if and only if condition {\rm (\ref{eq.CR})} is fulfilled and there is a 
section ${\mathcal F }_0\in H(\Omega, E_0) $, coinciding with $F_0$ in  
$D^+$ and such that $\Delta_0 {\mathcal F }_0 =0$ in $\Omega$. 
\end{corollary}

{\bf Proof.} If $i=0$ then the operator $A_{-1}^*$ in (\ref{eq.DeltaF}) 
equals to zero and therefore $\Delta_0 {\mathcal F}_0 =0$ in $\Omega$.

Back, as $\Delta_0 {\mathcal F}_0=0$ then the section ${\mathcal F}_0$ is smooth 
in $\Omega$. According to \cite[Theorem 9.4.8]{Tark37} the section  ${\mathcal F}_0$
belongs to $H(\Omega, E_0)$ if and only if it has finite order of growth  
near $\partial \Omega$. As $D\subset \Omega$, the section ${\mathcal F}^-_0$ 
has the same order of growth (in $D$) near $\partial D$. Then ${\mathcal 
F}^-_0 \in H (D, E_0 )$, $u = F^-_0 - {\mathcal F}^-_0$ in  
$H (D, E_0)$ and $(F_0 - {\mathcal F}_0)\in \chi_D (H (D, E_0)$ 
because $F_0 ={\mathcal F}_0$ in $D^+$.
\hfill $\square$

In the next section we will obtain a similar result in positive degrees 
of the complex $\{A_i\}$ over Lebesgue space $L^2 (D,E_i)$ choosing a 
canonical solution $u$ in (\ref{eq.F0}). In any case, Theorem \ref{t.AKy} can 
be easily reformulated to be like Corollary \ref{c.AKy.ell} 

\begin{corollary} \label{c.AKy.cohom} 
The Cauchy Problem {\rm \ref{pr.Cau}} is solvable if and only if condition 
{\rm (\ref{eq.CR})} is fulfilled and there is a section ${\mathcal F }_i\in 
H(\Omega, E_i)$ such that $(F_i - {\mathcal F}_i) \in  \chi_D (H(D, E_i )) $ 
and  $\Delta_i {\mathcal F }_i $ co-homological to zero in $\Omega$ with 
respect to the complex $\{A_i\}$.
\end{corollary}

{\bf Proof.} It follows from Theorem \ref{t.AKy} and (\ref{eq.DeltaF}) because 
$A_i \circ A_{i-1} \equiv 0$.
\hfill $\square$

\section{The Cauchy problem in the Lebesgue space} \label{s.L2}

Consider now the case $s=0$. Denote $\Sigma_0$ the null space of the Cauchy 
Problem {\rm \ref{pr.Cau}} for $s=0$, i.e. $\Sigma_0$ consists of  $L^2(D, 
E_i)$-sections $w$ with $A_i w=0$ in $D$ and $\tau_i (w)=0$ on $\Gamma$, or, 
the same  
\begin{equation} \label{eq.Cau.0}
(w, A_i ^*v)_D = 0 \mbox{ for all } v \in C^\infty _{comp}(D\cup \Gamma, 
E_{i+1}).
\end{equation}
Formula (\ref{eq.Cau.0}) guarantees that $\Sigma_0$ is a (closed) subspace 
in $L^2(D, E_i)$.

As the adjoint complex $\{A_i^*\}$ is elliptic too we may give similar 
definition of weak boundary value of a normal part (with respect to 
$\{A_i\}$) of a section on $\Gamma$.

\begin{definition} \label{def.BV.nu} 
We say that a section $u \in H_{A^*_{i-1}} (D, E_i)$, satisfying $A_{i-1}^* u 
=h$ in $D$ with $h \in H (D, E_{i-1})$,  has a weak boundary value $\nu
_{i,\Gamma} (u) = \nu_i (u_0)$ on $\Gamma$ for $u_0\in  {\mathcal D}'
(\Gamma,E_i)$  if 
$$ 
(h,g)_D - (u,A_{i-1} g)_D =  \langle \star \tilde \tau_{i-1} (g), \nu_i(u_0)
\rangle_{\Gamma} \mbox{ for all }  g \in   C^\infty _{comp}(D\cup \Gamma, 
E_{i-1}) .
$$
\end{definition}

\begin{theorem} \label{t.regular.sol0}
Let $f \in H^{-1} (D, E_{i+1})$, $u^0 =0$. If the Cauchy Problem {\rm \ref
{pr.Cau}} is solvable in $H^0_{A_i} (D, E_i )$ then its unique $L^2(D,E_i)$
-orthogonal to  $\Sigma_0$ solution $u (f)$ satisfies $\nu _{i,\Gamma}(u(f))=0$ 
on $\Gamma$ in the sense of Definition \ref{def.BV.nu} and $A_{i-1}^* u (f) =
0$  in $D$. 
\end{theorem}

{\bf Proof.} Obviously,  $H^0_{A_i} (D, E_i)$ is a Hilbert space with the 
scalar product 
$$
(\cdot ,\cdot )_{0,A_i}=(\cdot ,\cdot )_{0} + (A_i \cdot, A_i \cdot )_{-1}.
$$ 
Then the orthogonal complement to $\Sigma_0$ in this space coincides with 
$L^2(D,E_i)$-orthogonal complement to $\Sigma_0$. Thus, if the Cauchy Problem 
\ref{pr.Cau} has a solution $u$ in $H^0_{A_i}(D, E_i)$ then its $L^2(D,E_i)$
-orthogonal projection $u (f)$ to the orthogonal complement to $\Sigma_0$ is 
also a solution to Problem {\rm \ref{pr.Cau}} (it is evidently unique with the 
prescribed property). Clearly, any section of the type $A_{i-1} \phi$, with $v 
\in C^\infty_{comp} (D,E_{i-1})$, belongs to $\Sigma_0$. Hence  
$$
(u(f), A_{i-1} v)_D=0 \mbox{ for all } v \in C^\infty _{comp}(D, E_{i-1}),
$$ 
and then  $A_{i-1}^* u(f)=0$ in $D$. 

Now, according to Corollaries \ref{c.equiv.strong.nui} and \ref
{c.strong-weak.equiv.nu}, the section $u(f)$ has traces of $\nu_i (u (f))$ 
on $\partial D$, belonging to $H^{-1/2}(D,E_i)$. Hence, by Definition  
\ref{def.BV.nu}, the normal part $\nu_{i,\Gamma} (u)f))$ vanishes on 
$\Gamma$ if and only if 
\begin{equation} \label{eq.v.nu} 
(u(f), A_{i-1} v)_D  =0 \mbox{ for  all } v \in C^\infty _{comp}(D \cup 
\Gamma, E_{i-1}). 
\end{equation}
Further, it follows from Corollary \ref{c.equiv.strong.nui} that the space  
$H^0_{A_i} (D, E_i)$ is the Hilbert space with the scalar product 
$$
(\cdot , \cdot )_{0,\tau_i} = (\cdot , \cdot )_{0} + (\tau_i \cdot,\tau_i 
\cdot )_{-1/2}.
$$
Again we see that the orthogonal complement to $\Sigma_0$ in this space 
coincides with $L^2(D,E_i)$-orthogonal complement to $\Sigma_0$.
Denote  $\pi_{\tau_\Gamma}$ the orthogonal projection on the subspace 
$\Sigma_{\tau_{\Gamma}}$, consisting of sections with vanishing tangential 
parts on $\Gamma$. Definition \ref{def.BV.tau} guarantees that the subspace 
$\Sigma_{\tau_{\Gamma}}$ is closed in $H^0_{A_i} (D, E_i )$. As  $\tau   
_{i,\Gamma}( u(f))= u^0=0$ then for all $v \in C^\infty _{comp} (D \cup 
\Gamma, E_{i-1})$ we obtain:  
\begin{equation} \label{eq.v.nu2}  
(u(f), A_{i-1} v)_D = (\pi _{\tau_{\Gamma}} u(f), A_{i-1} v )_{0,\tau_i}  = 
( u(f), \pi _{\tau _{\Gamma}} A_{i-1} v)_{0,\tau _i} = ( u(f), \pi 
_{\tau_{\Gamma}} A_{i-1} v )_{D}. 
\end{equation} 
On the other hand, for all $g \in C^\infty _{comp}(D, E_{i+1})$ we have: 
$$
(\pi _{\tau _{\Gamma}} A_{i-1} v, A_i ^* g)_{D} = (\pi _{\tau _{\Gamma}} 
A_{i-1} v, A_i  ^* g)_{0,\tau_i}= (  A_{i-1} v, \pi _{\tau_{\Gamma}} A_i ^* 
g)_{0, \tau_i} = (A_{i-1} v, A_i ^* g)_{D} =0,
$$
because $A_{i} \circ A_{i-1}\equiv 0$. Therefore $A_i \pi _{\tau_{\Gamma}} 
A_{i-1} v=0$ in  $D$, and  $\pi _{\tau_{\Gamma}} A_{i-1} v \in \Sigma_0$ for 
all $v \in C^\infty _{comp}(D \cup \Gamma, E_{i-1})$. Hence, formulae (\ref
{eq.v.nu}) and (\ref{eq.v.nu2}) and the fact that $u(f)$ is orthogonal to 
$\Sigma_0$ in $L^2(D, E_i)$, imply that $\nu _{i,\Gamma}(u(f))=0$ on $\Gamma$.
\hfill $\square$


\begin{corollary} \label{c.regular.sol}
Let $f \in H^{-1} (D, E_{i+1})$, $u^0 \in H^{-1/2}(\overline \Gamma, E_i)$. If 
the Cauchy Problem {\rm \ref{pr.Cau}} is solvable in the space  $H^0_{A_i} 
(D, E_i)$ then the section $u (f,\tilde u_0) = u (f - A_i {\mathcal P}
_{\Delta_i} ^{(D)} \tau_i(\tilde u^0)) + {\mathcal P}_{\Delta_i}^{(D)}
\tau_i(\tilde u^0)$ is also its solution satisfying $\nu _i (u (f, \tilde 
u_0))=0$ on $\Gamma$, $ A_i ^* u (f,\tilde u_0) = A_i^*  {\mathcal P}_
{\Delta_i}^{(D)} \tau_i(\tilde u^0)$ in $D$. Besides, if $f \in H^{s} _{loc}(D 
\cup  \Gamma, E_{i+1})$, $u^0 \in H^{s+1/2}_{loc} (\Gamma, E_i)$ then $u (f, 
\tilde u_0) \in H^{s+1}_{loc} (D \cup \Gamma, E_i)$, $s \in {\mathbb Z}_+$.
\end{corollary}

{\bf Proof.} Let $u\in H^0_{A_i} (D, E_i )$ be a solution to Problem 
\ref{pr.Cau} with data $f \in H^{-1} (D, E_{i+1})$, $u^0 \in 
H^{-1/2}(\overline \Gamma, E_i)$. Then, according to Lemma \ref{l.projections} 
and Remark \ref{r.Dirichlet}, we have on $\Gamma$: 
$$
\tau _i (\tilde u^0) = \tau_i (u^0), \quad \tau_i ({\mathcal P}_{\Delta_i}
^{(D)} \tau_i (\tilde u^0)) =  \tau_i (u^0), \quad \nu _i ({\mathcal P}
_{\Delta_i}^{(D)} \tau_i (\tilde u^0)) =  0.
$$  
Hence Problem \ref{pr.Cau} with data $\hat f = f - A_i {\mathcal P}_{\Delta_i}
^{(D)}\tau_i (\tilde u^0) \in H^{-1} (D, E_{i+1})$ and $\hat u^0= 0$ is 
solvable in the space $H^0_{A_i} (D,E_i)$; the section $\hat u = u - {\mathcal 
P}_{\Delta_i}^{(D)} \tau_i(\tilde u^0)$ is its solution. Therefore Theorem \ref
{t.regular.sol0} implies that the section $u(f, \tilde u^0)$ is a solution to 
Problem \ref{pr.Cau} with data $f \in H^{-1} (D,E_{i+1})$, $u^0 \in 
H^{-1/2}(\overline \Gamma, E_i)$. By the construction it satisfies $\nu_i (u
(f, \tilde u_0))=0$ on $\Gamma$, $ A_i^* v (f, \tilde u_0) = A_i^*  {\mathcal 
P}_{\Delta_i}^{(D)} \tau_i (\tilde u^0)$ in $D$.

Finally, if $f \in H^{s} _{loc}(D \cup  \Gamma, E_i)$ then, using Theorem \ref 
{t.regular.sol0} and Lemma \ref{l.projections}, we conclude that $t(u(f))=0$ 
on $\Gamma$, 
$$ 
\Delta_i u(f) = (A_i ^* A_i + A_{i-1} A_{i-1}^*) u(f) = A_i ^* f  
\in H^{s-1} _{loc}(D \cup  \Gamma, E_i). 
$$
Therefore $u (f) \in H^{s+1}_{loc} (D \cup \Gamma, E_i)$, $s\in {\mathbb Z}_
+$, because of Theorem on local improvement of smoothness for solutions to 
Dirichlet Problem (see, for instance, \cite[Theorem 9.3.17]{Tark37}). 
Similarly, if $u^0 \in H^{s+1/2}_{loc} (\Gamma, E_i)$ then ${\mathcal 
P}_{\Delta_i}^{(D)} \tau_i( \tilde u^0) \in H^{s+1}_{loc} (D \cup \Gamma, E_i)$ 
according to Remark \ref{r.Dirichlet} and \cite[Theorem 9.3.17]{Tark37}). 
Thus, $u(f, \tilde u_0)$ belongs to  $H^{s+1}_{loc} (D \cup \Gamma, E_i)$, 
which was to be proved.
\hfill $\square$

Since Corollary \ref{c.regular.sol} practically reduces the Cauchy Problem 
\ref{pr.Cau} to the case with zero boundary data, we consider the situation in 
detail.

\begin{theorem} \label{t.AKy0}
Let $\Delta_{i-1}$, $\Delta_i$, $\Delta_{i+1}$ satisfy the Uniqueness 
Condition {\rm \ref{US}}. If $f \in H^{-1} (D, E_{i+1})$, $u^0 =0 $ then 
Problem {\rm \ref{pr.Cau}} is solvable in the space $H^0_{A_i} (D, E_i)$ if 
and only if $A_{i+1} f = 0$ in $D$, $\tau_{i+1}(f) = 0$ on $\Gamma$ and there 
is a section ${\mathcal F}_i\in L^2 (\Omega, E_i)\cap S_{\Delta_i} (\Omega) $ 
coinciding with $T_if$ in $D^+$.
\end{theorem}

{\bf Proof.} As $u^0=0$, then $F_i= T_if$. Moreover, by Lemma \ref{l.CR},  
condition (\ref{eq.CR}) is equivalent to the following two conditions: 
$A_{i+1} f = 0$ in $D$ and $\tau_{i+1} (f) = 0$ on $\Gamma$. Now if there is a 
section ${\mathcal F}_i \in L^2 (\Omega, E_i)\cap S_{\Delta_i} (\Omega)$ 
coinciding with $T_if$ in $D^+$ then $(T_if)^\pm, {\mathcal F}_i^\pm \in L^2 
(D^\pm , E_i)\cap S_{\Delta_i}(\Omega)$, $(T_if-{\mathcal F}_i)\in\chi_D
(L^2(D, E_i))$ and $A_i \Delta_i {\mathcal F }_i =0$ in $\Omega$. Therefore, 
it follows from Corollary \ref{c.AKy} that Problem \ref{pr.Cau} is solvable in 
the space $H^0_{A_i} (D, E_i)$ if $A_{i+1} f = 0$ in $D$ and $\tau _{i+1} (f) 
= 0$ on $\Gamma$. We note that formulae (\ref{eq.F0}) and (\ref{eq.sol}) yield:
\begin{equation} \label{eq.sol0}
u(f) = (T_{i} f - {\mathcal F}_i^- )\in L^2 (D, E_i).
\end{equation}  

Back, if Problem \ref{pr.Cau} is solvable in the space $H^0_{A_i} (D, E_i)$ 
then $A_{i+1} f = 0$ in $D$, $\tau_{i+1} (f) = 0$ on $\Gamma$. Moreover, 
the extension ${\mathcal F}_{i,u} \in L^2 (D, E_i ) \cap S_{A_i \Delta_i } 
(\Omega)$ of the section $T_if$ from  $D^+$ on  $\Omega$ is given by  
formula (\ref{eq.F0}). Putting the solution $u(f)$ into (\ref{eq.F0}) and 
using formula (\ref{eq.DeltaF}) and Definition \ref{def.BV.nu}, we obtain for 
all $v \in C^\infty _{comp}(\Omega, E_i)$: 
$$
- \langle \Delta_i {\mathcal F}_{i,u(f)}, v \rangle _\Omega = \langle \chi_D 
u(f) , A_{i-1} A_{i-1} ^*v  \rangle _\Omega = (u(f), A_{i-1} A_{i-1} ^* v)_D =
(A_i^*  u(f), A_i ^* v)_D =0 ,
$$
because $\nu _{i,\Gamma}(u (f)) =0$,  $A_{i-1}^* u(f)=0$ in $D$, i.e. 
$\Delta_i  {\mathcal F}_{i,u(f)} =0$ in $\Omega$.
\hfill $\square$

\begin{remark} \label{r.sol}
Theorem \ref{t.AKy0} easily implies conditions of {\it local} solvability of 
the Cauchy problem for complex $\{A_i\}$ in $L^2 (D,E_i)$ for $u^0=0$. Indeed, 
fix a point $x_0 \in \Gamma$. Let $U$ be a (one-sided) neighborhood of $x_0$ 
in $D$ and $\hat \Gamma=\partial U \cap \Gamma$. Set $ \hat F_i = T_{i} \chi_U 
f $. As $ F_i= \hat F_i +   T_{i} \chi_{D\setminus U}f$ we see that $ F^+ _i$ 
extends as a solution to the Laplacian $\Delta_i $ in $\hat \Omega= U\cup 
\hat\Gamma \cup D^+$ if and only if the potential $\hat F^+_i $ does. Hence, 
under condition {\rm (\ref{eq.CR})}, the solution of the Cauchy problem exists 
in the neighborhood $U$ where the extension of the potential $F^+_i$ does.  
\end{remark}

Also we would like to note that Theorem \ref{t.AKy0} gives not only the
sol\-va\-bi\-li\-ty conditions to Problem \ref{pr.Cau} but the solution
itself, of course, if it exists (see (\ref{eq.sol0})). It is clear that we can 
use the theory of functional series (Taylor series, Laurent series, etc.) in 
order to get information about extendability of the potential $T_i^+f $ (cf. 
\cite{AKy}, \cite{Tark37}). However in this paper we will use the theory of 
Fourier series with respect to the bases with the double orthogonality 
property (cf. \cite{Shp}, \cite{Tark37} or elsewhere). Moreover, using formula 
(\ref{eq.sol0}) we can construct approximate solutions of problem \ref{pr.Cau} 
in the Lebesgue space $L^2 (D,E_i)$. 

\begin{lemma} \label{l.bdo}
If $\omega \Subset \Omega$
is a domain with a piece-wise smooth boundary and $\Omega \setminus
\omega$ has no compact (connected) components then there exists an orthonormal
basis $\{b_\nu \}_{\nu=1}^\infty$ in $L^2 (\Omega,E_i) \cap S_{\Delta_i} 
(\Omega)$ such that $\{b_{\nu|\omega}\}_{\nu=1}^\infty$ is an orthogonal basis 
in $L^2 (\omega,E_i) \cap S_{\Delta_i} (\omega)$.
\end{lemma}

{\bf Proof.} In fact, these $\{b_\nu \}_{\nu=1}^\infty$ are
eigen-functions of compact self-adjoint  linear operator $R(\Omega,
\omega)^*R(\Omega, \omega)$, where 
$$
R(\Omega, \omega): L^2 (\Omega,E_i) \cap S_{\Delta_i} (\Omega)
\to L^2 (\omega,E_i) \cap S_{\Delta_i} (\omega)
$$ is the natural inclusion operator (see
\cite{Tark37} or \cite[theorem 3.1]{ShTaLMS}). \hfill $\square$

Now we can use the basis $\{ b_\nu \}$ in order to simplify Theorem \ref
{t.AKy0}. For this purpose fix domains $\omega \Subset D^+$ and $\Omega$ as in 
Lemma \ref{l.bdo} and denote by 
$$
c_\nu (T_i f^+) = \frac{(T_i f^+, b_\nu )_{L^{2} (\omega,E_i)}}{\|b_\nu
\|^2_{L^{2}(\omega,E_i)}}, \quad \nu \in {\mathbb N},
$$ 
the Fourier coefficients of $T_i f^+$ with respect to the
orthogonal system $\{ b_{\nu|\omega}\}$ in $L^2 (\omega,E_i)$.

\begin{corollary} \label{c.bdo}
Let $f \in H^{-1} (D, E_{i+1})$, $u^0 =0 $. Problem {\rm \ref{pr.Cau}} is 
solvable in the space $H^0_{A_i} (D, E_i)$ if and only if $A_{i+1} f = 0$ in 
$D$, $\tau_{i+1}(f) = 0$ on $\Gamma$ and and the series $\sum_{\nu = 1}^\infty 
|c_\nu (T_if^+)|^2$ converges.
\end{corollary}

{\bf Proof.} Indeed, if Problem \ref{pr.Cau} is solvable in $L^2 (D,E_i)$
then, according to Theorem  \ref{t.AKy0} condition {\rm (\ref{eq.CR})} is 
fulfilled, and there  exists a function ${\mathcal F}_i \in L^2 (\Omega,E_i) 
\cap S_{\Delta_i} (\Omega)$ coinciding with $T_i f^+$ in $\omega$.
By Lemma \ref{l.bdo} we conclude that 
\begin{equation} \label{eq.bdo} 
{\mathcal F}_i (x) = \sum_{\nu = 1}^\infty k_\nu (
{\mathcal F}_i) b_\nu (x), \quad x \in \Omega, 
\end{equation}
where $k_\nu ({\mathcal F}_i) = ({\mathcal F}_i, b_\nu )_{L^2 (\Omega,E_i)}$, 
$\nu \in {\mathbb N}$, are the Fourier coefficients of ${\mathcal F}_i$ with 
respect to the orthonormal basis $\{ b_{\nu}\}$ in $L^2 (\Omega,E_i) \cap 
S_{\Delta_i} (\Omega)$. Now Bessel's inequality implies that the series
$\sum_{\nu = 1}^\infty |k_\nu ({\mathcal F}_i)|^2$ converges.

Finally, the necessity of the corollary holds true because 
$$
c_\nu (T_i f ^+) =\frac{(R(\Omega, \omega){\mathcal F}_i, R(\Omega,
\omega)b_\nu )_{L^2 (\omega,E_i)}}{(R(\Omega, \omega)b_\nu , R(\Omega, 
\omega)b_\nu )_{L^2(\omega,E_i)}} = \frac{({\mathcal F}_i, R(\Omega, 
\omega)^*R(\Omega,\omega)b_\nu )_{L^2 (\Omega,E_i)}}{(b_\nu , R(\Omega, 
\Omega)^*R(\Omega, \omega)b_\nu )_{L^2(\omega,E_i)}} = k_\nu ({\mathcal F}_i).
$$

Back, if the hypothesis of the corollary holds true then we invoke the 
Riesz-Fisher theorem. According to it, in the space $L^2 (\Omega,E_i) 
\cap S_{\Delta_i} (\Omega)$ there is a section
\begin{equation} \label{eq.Phi}
{\mathcal F}_i (x) = \sum_{\nu = 1}^\infty c_\nu (T_if ^+) b_\nu (x), \qquad x 
\in \Omega.
\end{equation}
By the construction, it coincides with $T_i f ^+$ in $\omega$. Therefore, 
using Theorem \ref{t.AKy0}, we conclude that Problem \ref{pr.Cau} is solvable 
in $L^2(D,E_i)$.
\hfill $\square$

The examples of bases with the double orthogonality property be found in
\cite{ShTaLMS}, \cite{Tark37}, \cite{Shp}.

Let us obtain Carleman's formula for the solution of Problem \ref{pr.Cau}.
For this purpose we introduce the following Carleman's kernels:
$$
{\mathfrak C}_N (y ,x ) = (A_i ^*)^\prime _y \Phi_i (y, x) - \sum_{\nu=1}^N 
c_\nu ((A_i ^*)^\prime _y \Phi_i (y , \cdot)) b_\nu (x),  \, N \in {\mathbb 
N}, \, x \in \Omega, \,y \not \in \overline \omega, x \ne y.
$$

\begin{corollary} \label{c.sol}
If Problem  \ref{pr.Cau} is solvable in $L^2 (D,E_i)$ for data $u_0 =0$ 
and $f \in L^2 (D, E_{i+1}) \cap H^s _{loc} (D \cup \Gamma, E_{i+1})$ then $u 
(f)$ belongs to $ H^{s+1} _{loc} (D \cup 
\Gamma, E_{i})$ and the following Carleman formula holds:
\begin{equation} \label{eq.Carl.1}
u(f) (x)=  \lim_{N\to \infty} \int_{D}   \langle {\mathfrak C}_N 
(\cdot , x ) , f \rangle _y \ dy 
\end{equation}
where the limit converges in the spaces $ H^0_{A_i} (D,E_i)$ and $H^{s+1}_{loc} (D \cup \Gamma, E_{i})$.
\end{corollary}

{\bf Proof.} Since $\overline \omega \cap \overline D = \emptyset$, 
using Fubini Theorem we have for all $\nu \in \mathbb N$:
\begin{equation*} 
c_\nu(T_if^+) = \int_{D}   \langle c_\nu((A_i ^*)^\prime _y \Phi_i (y, \cdot 
)) , f  \rangle _y \ dy .
\end{equation*}
This exactly yields identity (\ref{eq.Carl.1}) after applying Corollary \ref
{c.bdo},  formula (\ref {eq.Phi}) and regrouping the summands in 
(\ref{eq.sol0}). 

Besides, since ${\mathcal F}_i$ and each function $b_\nu$ are solutions of the 
elliptic system $\Delta_i$ in $\Omega$, the Stiltjes-Vitali theorem implies 
that the series (\ref{eq.Phi}) converges in $C^\infty_{loc} (\Omega,E_i)$.  
Therefore we additionally conclude that the limit converges to $ u (f)$ in 
$H^{s+1}_{loc} (D \cup \Gamma,E_i)$ because $T_if \in H^{1} (D,E_i) \cap 
H^{s+1} _{loc} (D \cup \Gamma, E_{i}) $ due to the transmission property (see 
\cite{ReSch}).
\hfill $\square$ 

Considering general complexes with smooth coefficients we arrive to the following 
natural question: under what conditions on the domain $D$ the complex  $\{A_i\}$
is exact at the positive degrees ? As far as we know there is no answer in 
the general situation. It is known that the formally exact differential  elliptic complexes 
with real analytic coefficients are locally exact at the positive degrees 
(see, for instance, \cite{Ta1}, \cite{Spe}). Of course, all the Hilbert complexes 
with constant coefficients are exact at the positive degrees over the spaces of distrubutions 
in convex domains (see, for instance, \cite{Pal}. Thus we are to consider this 
most investigated situation. However we emphasize that the use of the above proposed 
approach to the Cauchy problem for the elliptic complexes does not involve the information 
on the exactness of the complex!

\section{Complexes with constant coefficients}

Now we are to discuss examples for complexes with constant coefficients. 
Actually we can say much more, at least for domains of the special type.

\begin{corollary} \label{c.Carleman} 
Let (\ref{eq.complex}) be an elliptic first order complex with constant 
coefficients in ${\mathbb R}^n$. If $\partial D \setminus \Gamma$ is a part of 
a strictly convex domain $\Omega \supset D$, then for any section $w \in  
C^\infty (\overline D, E_i)$ there is a section $h\in L^2 (D, E_{i-1})\cap  
C^{\infty}_{loc} (D\cup \Gamma, E_{i-1})$ such that $\tau_i (A_{i-1} h) =0$ on 
$\Gamma$ and the following formula holds true:
\begin{equation} \label{eq.Carl.2}
w(x)=  {\mathcal P}^{(D)}_{\Delta_i} \chi_{\Gamma}\tau_i(w) )x) +  \lim_{N\to 
\infty} \int_{D} \langle {\mathfrak C}_N (\cdot , x ) ,  A_i (u -  {\mathcal P}^{(D)}_{\Delta_i} \chi_{\Gamma}\tau_i(w)) \rangle_y \ dy  +  A_{i-1} h (x), 
\end{equation}
where the limit converges in the spaces $H^0_{A_i} (D,E_i)$ and $C^\infty 
_{loc}(D\cup \Gamma, E_i)$.
\end{corollary}

{\bf Proof.} Under the hypothesis of the corollary, Problem \ref{pr.Cau} is 
solvable for the data  $\tau_{i,\Gamma}(u) \in C^{\infty} (\overline \Gamma, E_i)$ 
and $A_i u \in C^\infty (\overline D, E_{i+1})$. Extending $\tau_{i,\Gamma}(w)$
by zero onto all the boundary of $D$, we obtain $\tilde w_0 =  \chi_{\Gamma}
\tau_i(w) \in L^2 (\partial D, E_i)$. Now Corollary \ref{c.regular.sol} 
implies that the section $u ( A_i w, \tilde w_0)$ belongs to the space 
$L^2 (D, E_i)\cap  C^{\infty}_{loc} (D\cup \Gamma,E_i)$ and $v= w- u (A_i w, 
\tilde u_0) \in \Sigma _0 \cap  C^{\infty}_{loc} (D\cup \Gamma,E_i)$. 

Denote $v_0$ the extension by zero of $v$ from  $D$ on $\Omega$. Clearly, 
$v_0\in L^2 (\Omega, E_i)$. As $\tau _{i,\Gamma}(v)=0$, then  
$A_{i+1} v_0 = 0$ in $\Omega$ and hence there is a section $\tilde h 
\in L^2 (\Omega, E_{i-1})\cap H^1_{loc} (\Omega,E_{i-1})$ such that  
$A_{i-1} \tilde h=v_0$ in $\Omega$ (see, for instance, \cite{Pal}). 
Set $h= \tilde h - A_{i-2} \Phi _{i-2}  \chi_D A_{i-2}^* \tilde h $.
Then 
$$ 
A_{i-2}^* \Delta _{i-1} h =  A_{i-2}^* A_{i-2} A_{i-2}^* 
\tilde h - A_{i-2}^* A_{i-2} \chi_D A_{i-2}^* \tilde h = 0 \mbox{ in } D,  
$$
$$ 
A_{i-1} \Delta _{i-1} h =  A_{i-1} A_{i-1} ^*  A_{i-1}   
\tilde h  = A_{i-1} A_{i-1} ^*  v  \mbox{ in } D,  
$$
As the operator $(A_{i-1}  \oplus A_{i-2}^*) \Delta _{i-1}$ has injective 
symbol and 
$$
(A_{i-1}  \oplus A_{i-2}^*) \Delta _{i-1} h = ( A_{i-1} A_{i-1} 
^*  v ,0) \in C^\infty _{loc} (D \cup \Gamma, (E_{i}, E_{i-2})), 
$$ 
we see that $h \in C^\infty _{loc} (D \cup \Gamma, E_{i-1})$ satisfies 
$A_{i-1} h=w$  in $D$. Thus, $v = u (A_i w, \tilde u_0) + A_{i-1} h$ and 
formula (\ref
{eq.Carl.2}) follows from Corollary \ref{c.sol}. 
\hfill $\square$

At the conclusion let us consider two examples. 

\begin{example} \label{ex.d}
{\rm
Let (\ref{eq.complex}) be the de Rham complex over ${\mathbb R}^n$, i.e $E_i$ be 
the bundle of the exterior differential forms of the degree $i$ and $A_i$ be 
the differentiation operator $d_i$ for the exterior  
differential forms. Choosing coordinates $x=(x_1,...,x_n)\in{\mathbb R}^{n}$ we 
have for a form $u \in C^\infty ({\mathbb R}^n, \Lambda^{i})$:
$$
u = \sum_{|I|=i}  u_I (x) \ d x_I, \quad 
d_i u = \sum_{j=1}^n \sum_{|I|=i} \frac{\partial u_I  }{ 
\partial x_j} (x)\ d x_j \wedge d  x_I,
$$ 
where  $I = (j_1, \dots, j_i)$, 
$d x_I = d  x _{j_1} \wedge \dots \wedge x_{j_i} $  and $\wedge$ is the exterior product 
for the differential forms. 

Let  $*$ be the Hodge operator for the differential forms (see, for instance, 
\cite{Ta1}), in particular, $dx_I \wedge *dx_I =dx $.  
Then $\Delta_i = \Delta I_{k(i)}$, where $\Delta$ is
the usual Laplace operator in ${\mathbb R}^{n}$ and $I_{k(i)}$ is the unit 
$k(i) \times k(i)$-matrix. If $\Phi_i = I_{k(i)} \ \Phi$, where $\Phi$ is the 
standard fundamental solution to $\Delta$ of the convolution type, then $M_i$ 
is the Norguet integral and  (\ref{eq.MB}) is the the 
Norguet integral formula (see, for instance,  \cite[\S 2.5]{Ta1}). 

Let $\{ h^{(j)}_\nu \}$  be the set of homogeneous harmonic polynomials 
forming a complete orthonormal system in the space $L^2 (\partial 
B(0,1) )$ on the unit sphere $\partial B(0,1)$ in ${\mathbb R}^{n}$, $n\geq 2$ 
(see \cite[p. 453]{So}). Therefore $\{ h^{(j)}_{\nu |\partial B(0,1)} \}$ are 
{\it spherical harmonics} where $\nu$ is the homogeneity, $j$ is the number of 
the polynomial of degree $\nu$ in the basis, $1\leq j \leq J(\nu ,n)$ with  
$J(\nu ,n)=\frac{(n+2\nu -2)(n+\nu -3)!}{\nu !(n-2)!}$, $\nu>0$, $J (0,n)=1$. 
It is easy to see that the system $\{ h^{(j)}_\nu \}$ is orthogonal in  
$L^2 (B(0,R))$ for any ball $B(0,R)$.

Let $D$ be a part of the unit ball $\Omega$ cut off by a hypersurface $\Gamma 
\not \ni 0$. Then Carleman kernel in formulae (\ref{eq.Carl.1}), (\ref
{eq.Carl.2}) has the following form:
$$
{\mathfrak C}_N (y, x) =   * _{y} d _{y} 
\Phi_i (y, x)- \sum_{|I|=i} \sum_{\mu =0}^N \sum_{j =1}^{J(\mu,n)} * 
_{y} d_{y} \left( \frac{{h_{\mu} 
^{(j)}(y)}  \ d y _I  }{|y|^{n+2\mu -2}(n+2\mu -2)} 
\right) \ h_{\mu} ^{(j)}(x) \ d x_I. 
$$  

We note that the operators $d_i$ are non-zero for $0 \leq i \leq n-1$ only.

Hence for $n=1$ the operator $d_0$ is the usual differentiation and all the other operators 
$d_i$ are identically zeros. Then the Cauchy problem for an interval $D = (a,b) \subset 
{\mathbb R}$ is well known: given a distribution $f$ on $(a,b)$ find a distribution 
$u$ on $[a,b)$ such that 
$$
\left\{ 
\begin{array}{lll}  u' (x) = f (x), & x \in (a,b),  \\
u(a) = 0.\\ 
\end{array}
\right. 
$$  
This problem is well-posed in the Sobolev spaces and its solution is given by the integral 
$$
u (x) = \int_a^x f(t) \ dt,
$$
at least for $f$ from the Sobolev spaces of a non-negative smoothness. For 
elements $f$ from the Sobolev spaces of a negative smoothness the interpretations 
of the integral are also well known.

For $n=2$ the Cauchy problem for the de Rham complex at the degree $i=1$ can be 
inerpretated as follows. Let $D$ be a bounded domain in ${\mathbb R}^2 $ and  
$$
G = \{(x_1,x_2, x_3): (x_1,x_2) \in D, \, 0<x_3 <A\} \subset {\mathbb R}^3
$$ 
be a cylinder with the base $D$. If we consder $G$ as a bassin where the  
liquid behaves similarly in every section 
$$D_b = \{(x_1,x_2, b): (x_1,x_2) \in D, \, 0<b <A\} 
$$
then the (stationary) flow of the ideal non-contractible liquid can be 
described by the system of equations 
$$
\left\{
\begin{array}{ccc} 
 \frac{\partial u_1}{\partial x_2} - \frac{\partial u _2}{\partial x_1} = h & \mbox{ in } & D, \\
\frac{\partial u_1}{\partial x_1}  + \frac{\partial u_2}{\partial x_2} = g & \mbox{ in } & D,\\
\end{array}
\right.
$$
where the vector $u = (u_1,u_2)$ corresponds to the velocity  vector of the fluid and 
the components $h$, $g$ reflect the rotation points and the source points 
respectively (see, for example, \cite[Ch. III, \S 2]{LaShb}). This exactly means 
$$
d_1 u = f  \mbox{ in }  D, \, d_0^* u = - g \mbox{ in } D,
$$  
for the differential forms 
$$
u (x)= u_1 (x) dx_1 + u_2 (x) dx_2, \quad f(x) = h (x) dx_1 \wedge dx_2 , \quad g (x)  
$$
of the degrees $1$, $2$ and $0$ respectively. If  $(n_1 (x), n_2 (x))$ is the unit 
normal vector with respect to $\partial D$ at the point $x$ then 
$$
\tau_1 (u) = n_2 u_1 - n_1 u_2 \mbox{ on }  \partial D, \quad \nu_1 
(u) = n_1 u_1 + n_2 u_2 
\mbox{ on }  \partial D . 
$$
According to  Theorem \ref{t.AKy0} the Cauchy problem for the de Rham complex 
in $D$ with boundary data on $\Gamma \subset \partial D$, i.e. 
$$
\left\{
\begin{array}{ccc} 
d_1 u = f & \mbox{ in } &  D, \\
\tau _1 (u)=0 & \mbox{ on }  & \Gamma, \\
\end{array}
\right.
$$
is equivalent to the following problem 
$$
\left\{
\begin{array}{ccc} 
d_1 v = f & \mbox{ in } &  D, \\
d^*_0 v = 0 & \mbox{ in } &  D, \\
\tau _1 (u)=0 & \mbox{ on }  & \Gamma, \\
\nu_1 (u) = 0 & \mbox{ on } &  \Gamma. \\
\end{array}
\right.
$$
The last one is obviously the Cauchy problem for the classical Cauchy-Riemann system 
with respect to the function $w(z) = v_2 (x_1,x_2) + \sqrt{-1} v_1 (x_1,x_2)$ with 
$z=x_1 + \sqrt{-1} x_2$:
$$
\left\{
\begin{array}{ccc} 
\frac{\partial w}{\partial \overline z} =  f/2 & \mbox{ in } &  D, \\
w = 0 & \mbox{ on } &  \Gamma, \\
\end{array}
\right.
$$
where $\frac{\partial }{\partial \overline z} = \frac{1}{2}(\frac{\partial }{\partial x_1} + 
\sqrt{-1}\frac{\partial }{\partial  x_2})$. 
Thus, according to Hadamard's example (see \cite{Hd}) the Cauchy problem for the de Rham complex 
in ${\mathbb R}^2$ at the degree $1$ is ill-posed in all the standard functional spaces (the 
spaces of smooth functions, the Sobolev spaces etc.). 

For $n=3 $ the operators $d_0$, $d_1$, $d_2$ 
can be identifyed with the famuos gradient operator $\nabla$, the rotor operator 
$\mbox{rot}$ and the divergence  operator $\mbox{div}$ respectively 
which are widely used in Mechanics, Hydrodinamics, Electrodynamics and so on:
$$ 
d_0 \approx 
\nabla = \left( 
\begin{array}{ccc} 
\frac{\partial}{\partial x_1}  \\
 \frac{\partial}{\partial x_2}  \\
 \frac{\partial}{\partial x_3}   \\
\end{array}
\right), \quad 
d_1 \approx  
\mbox{rot} = 
\left( 
\begin{array}{ccc} 
0 & - \frac{\partial}{\partial x_3} &  \frac{\partial}{\partial x_2} \\
 \frac{\partial}{\partial x_3} & 0 & - \frac{\partial}{\partial x_1} \\
- \frac{\partial}{\partial x_2} &  \frac{\partial}{\partial x_1} & 0  \\
\end{array}
\right), 
 \quad 
d_2\approx 
\mbox{div} = \left( 
\begin{array}{ccc} 
\frac{\partial}{\partial x_1}  &  \frac{\partial}{\partial x_2}  & 
 \frac{\partial}{\partial x_3}   \\
\end{array}
\right),
$$
$$
d^*_0\approx  - \mbox{div}, \quad d^*_2\approx \mbox{rot}, \quad d^*_2 \approx - \nabla
$$
For instance, according  to  Theorem \ref{t.AKy0}, 
the Cauchy problem for the de Rham complex at the degree $1$ 
for a domain $D \subset {\mathbb R}^n$, a set $\Gamma \subset \partial D$ and a datum  
$f=(f_1, f_2, f_3)$ is equivalent 
to the Cauchy problem for the (stationary) Maxwell type system  with respect to the 
vector function $u=(u_1, u_2, u_3)$:
$$
\left\{ 
\begin{array}{lll}   \mbox{rot} \, u = f   & \mbox{ in } & D ,\\
\mbox{div} \, u  = 0   & \mbox{ in } & D ,\\ 
u = 0 & \mbox{ on } & \Gamma.\\ 
\end{array}
\right.
$$  
We refer to \cite{MMT} for applications of the theory of differential complexes to 
the investigation of the Maxwell type equations.
}
\end{example}

\begin{example} \label{ex.bdo}
{\rm
Let (\ref{eq.complex}) be the Dolbeault complex over ${\mathbb C}^n$, i.e $E_i$ be 
the bundle of exterior differential forms of bi-degree $(0,i)$ and $A_i$ be 
the  Cauchy-Riemann operator $\overline \partial_i$ for the exterior  
differential forms. Choosing coordinates $z=(z_1,...,z_n)$ with $z_j=x_j+
\sqrt{-1}x_{j+n}$, $j=1,...,n$, and $x=(x_1,...,x_{2n})\in{\mathbb R}^{2n}$ we 
have for a form $u \in C^\infty ({\mathbb C}^n, \Lambda^{(0,i)})$:
$$
u = \sum_{|I|=i}  u_I (z) \ d \overline z_I, \quad 
\overline \partial_i u = \sum_{j=1}^n \sum_{|I|=i} \frac{\partial u_I  }{ 
\partial \overline z_j} (z)\ d \overline z_j \wedge d \overline z_I,
$$ 
where $ \frac{\partial}{\partial \overline z_j}  = \frac{1}{2} 
\left(\frac{\partial}{\partial x_j} +  \sqrt{-1} \frac{\partial}{\partial 
x_{j+n}}\right)$, $dz_j = dx_j + \sqrt{-1} dx_{j+n}$, $I = (j_1, \dots, j_i)$, 
$d \overline z_I = d  \overline z _{j_1} \wedge \dots \wedge  d  \overline z 
_{j_i} $. 

It is well known that $\star u =  \overline {* 
u}$ for a form $u$ with $*$ being the Hodge operator for the differential forms 
(see \cite[\S 14]{Ky}). 
Then $\Delta_i =1/2\ \Delta I_{k(i)}$, where $\Delta$ is
the usual Laplace operator in ${\mathbb R}^{2n}$ and $I_{k(i)}$ is the unit 
$k(i) \times k(i)$-matrix. If $\Phi_i = I_{k(i)} \ \Phi$, where $\Phi$ is the 
standard fundamental solution to $\Delta$ of the convolution type, then $M_i$ 
is the Martinelli-Bochner-Koppelmann integral and  (\ref{eq.MB}) is the the 
Martinelli-Bochner-Koppelmann integral formula (see, for instance, \cite{AYu} 
or \cite{Ta1}). 


Let $D$ be a part of the unit ball $\Omega$ cut off by a hypersurface $\Gamma 
\not \ni 0$. Then Carleman kernel in formulae (\ref{eq.Carl.1}), (\ref
{eq.Carl.2}) has the following form (see \cite{FeShl2}):
$$
{\mathfrak C}_N (\zeta, z) =   \star _{\zeta} \overline \partial _{\zeta} 
\Phi_i (\zeta, z)- \sum_{|I|=i} \sum_{\mu =0}^N \sum_{j =1}^{J(\mu,2n)} \star 
_{\zeta} \overline \partial _{\zeta} \left( \frac{\overline {h_{\mu} 
^{(j)}(\zeta)}  \ d\overline \zeta _I  }{|\zeta|^{2n+2\mu -2}(2n+2\mu -2)} 
\right) \ h_{\mu} ^{(j)}(z) \ d\overline z_I
$$  
where $\{h_\mu ^{(j)}\}$ is the system of the spherical harmonics (see Example \ref{ex.d}).

A result similar to Corollary \ref{c.Carleman} was obtained in \cite[Theorem 
3.1]{NaSchTa} for the Dolbeault complex if $\partial D\setminus \Gamma$ is 
$i$-strictly pseudo concave hypersurface; however they had no aim to prove 
that the tangential part of the rest $\overline \partial _i h$ vanished 
on $\Gamma$.    
}
\end{example}


\smallskip

\textit{Acknowledgments.\, The investigations  were supported 
by DAAD and by RFBR grant 11-01-00852a.}

\end{document}